\documentclass[11pt]{article}
\usepackage{amsmath}
\usepackage{amssymb}
\usepackage{amscd}
\usepackage{latexsym}
\usepackage{theorem}
\usepackage{doublespace}
\usepackage{epsfig}
\usepackage{psfrag}
\usepackage{amsfonts}
\usepackage{enumerate}
\usepackage{bbold}
\usepackage{comment}

\numberwithin{equation}{section}
\numberwithin{figure}{section}

\newtheorem{theorem}{Theorem}[section]
\newtheorem{lemma}[theorem]{Lemma}

\newtheorem{proposition}[theorem]{Proposition}

{\theorembodyfont{\rmfamily}
\theoremstyle{plain}
\newtheorem{definition}[theorem]{Definition}
\newtheorem{example}[theorem]{Example}

\newtheorem{remark}[theorem]{Remark}

}

\newenvironment{proof}{\noindent {\bf Proof:} 
\nopagebreak}{\hfill $\Box$ \vspace{2 ex}}

\usepackage{multicol}

\makeatletter

\def\theglossary{\@restonecoltrue\if@twocolumn\@restonecolfalse\fi
\columnseprule\z@ \columnsep 35\p@
\let\@makessectionhead\indexsec
\@xp\section\@xp*\@xp{\glossaryname}%
\let\item\@idxitem
\parindent\z@  \parskip\z@\@plus.3\p@\relax
\footnotesize}

\def\glossaryname{Notation Index}

\makeatother

\makeglossary

\input xy
 \xyoption{all}

\renewcommand{\to}{\rightarrow}
\newcommand{\into}{\hookrightarrow}

\newcommand{\Ker}{\operatorname{Ker}}

\newcommand{\Z}{{\mathbb Z}}

\newcommand{\Q}{{\mathbb Q}}
\newcommand{\R}{{\mathbb R}}
\newcommand{\C}{{\mathbb C}}

\newcommand{\g}{{\mathfrak g}}

\renewcommand{\t}{{\mathfrak t}}

\renewcommand{\mod}{{/\!\!/}}

\UseComputerModernTips

\DeclareMathOperator{\grad}{grad}

\DeclareMathOperator{\Stab}{Stab}
\DeclareMathOperator{\Spinc}{Spin^{c}}
\DeclareMathOperator{\Gr}{Gr}
\DeclareMathOperator{\Pic}{Pic}
\DeclareMathOperator{\BU}{BU}
\DeclareMathOperator{\Hom}{Hom}
\DeclareMathOperator{\End}{End}
\newcommand{\bbS}{\mathbb{S}}
\DeclareMathOperator{\rank}{rank}
\newcommand{\pt}{\mathrm{pt}}
\DeclareMathOperator{\Tor}{Tor}

\begin{document}
\begin{spacing}{1.1}

%Full addresss of first author:
%Megumi Harada 
%Department of Mathematics
%Bahen Centre for Information Technology
%University of Toronto 
%40 St. George St., Room 6290
%Toronto, Ontario M5S 2E4
%Canada 

\noindent
{\LARGE \bf Surjectivity for Hamiltonian $G$-spaces in $K$-theory}
\bigskip\\
{\bf Megumi Harada}\footnote{{\tt megumi@math.toronto.edu}}\\
Department of Mathematics, University of Toronto, 
Toronto, Ontario, M5S 2E4 Canada\smallskip \\
{\bf Gregory D. Landweber}\footnote{{\tt greg@math.uoregon.edu}
%\newline \mbox{~~~~} The second author would like to thank The Fields Institute and the University of Toronto for their hospitality.
\newline \mbox{~~~~} {\it MSC 2000 Subject Classification}: 
Primary: 53D20 \hspace{0.1in} Secondary: 19L47
% 53D20: Momentum maps; symplectic reduction
% 19L47: Equivariant K-theory
% 55R91: Equivariant fiber spaces and bundles 
% 57R22: Topology of vector bundles and fiber bundles
% DG, SG, KT 
\newline \mbox{~~~~} {\it Keywords:} equivariant $K$-theory, Kirwan surjectivity, Morse-Kirwan function, symplectic quotient, Atiyah-Bott lemma, equivariant formality}\\
Mathematics Department, University of Oregon,
Eugene, Oregon 97403-1222 USA\smallskip \\
\bigskip

%%%%%%%%%%%%%%%%%%%%
% Disclaimer
%%%%%%%%%%%%%%%%%%%%
%\begin{center}
%\framebox{
%{\Large\bf DRAFT (\today): DO NOT DISTRIBUTE.}}
%\end{center}

%%%%%%%%%%%%%%%%%%%%%
%  Abstract
%%%%%%%%%%%%%%%%%%%%%
{\small
\begin{quote}
\noindent {\em Abstract.}  Let $G$ be a compact connected Lie group,
and $(M,\omega)$ a Hamiltonian $G$-space with proper moment map $\mu$.  
We give a surjectivity
result which expresses the $K$-theory of the symplectic quotient $M
\mod G$ in terms of the equivariant $K$-theory of the original
manifold $M$, under certain technical conditions on $\mu$. This result
is a natural $K$-theoretic analogue of the Kirwan surjectivity theorem
in symplectic geometry. The main technical tool is the $K$-theoretic
Atiyah-Bott lemma, which plays a fundamental role in the symplectic
geometry of Hamiltonian $G$-spaces. We discuss this lemma in detail
and highlight the differences between the $K$-theory 
and rational cohomology versions
of this lemma.

We also introduce a $K$-theoretic version of equivariant formality and
prove that when the fundamental group of $G$ is torsion-free, every compact Hamiltonian $G$-space is equivariantly formal. Under
these conditions, the forgetful map $K_{G}^{*}(M)\to K^{*}(M)$ is surjective, and thus every complex vector bundle admits a stable equivariant structure. Furthermore, by considering complex line bundles, we show that every integral cohomology class in $H^{2}(M;\Z)$ admits an equivariant extension in $H_{G}^{2}(M;\Z)$.
%We also discuss equivariant formality in this setting, which in this
%$K$-theoretic context is the statement that every vector bundle can be
%given a $G$-equivariant structure. We prove a partial result for
%compact Hamiltonian $G$-spaces in the case where $\pi_1(G)$ has no
%torsion.
\end{quote}
}
\bigskip

\section{Introduction}

In this article, we explore the $K$-theoretic analogues of several
results of Kirwan from \cite{Kir84} involving Hamiltonian $G$-spaces,
where $G$ is a compact connected Lie group.  Our main result is
Theorem~\ref{thm:KthyKirwan}, 
which asserts that if $M$ is a Hamiltonian
$G$-space, then the $G$-equivariant $K$-theory of $M$ surjects onto
the ordinary $K$-theory of the symplectic quotient $M\mod G$. In other
words, the Kirwan map given by the composition
\[
\kappa : K_{G}^{*}(M) \to K_{G}^{*}(\mu^{-1}(0) ) \cong K^{*}(M\mod
G)
\]
is a surjection, provided that the moment map $\mu$ is proper and $0$ is a regular value. Kirwan's original theorem, which states that the analogous restriction
map $\kappa$ of rational cohomology rings is a surjection, is
crucial for the explicit description
of cohomology rings for many symplectic manifolds
(see for example \cite{Kir84}, \cite{Kir92}, and \cite[Chapter 8]{MFK94}).

%It is an open
%question whether results analogous to Kirwan surjectivity (as well as other
%related cohomological results) can be proven in 
%other contexts. For instance, we may consider other geometries possessing moment
%maps which give rise to interesting theories of quotients, such as
%contact geometry or hyperk\"ahler geometry. Similarly, we may consider 
%generalizations to other cohomology theories. The main result in this
%paper should therefore be considered as a first step toward a
%generalization of Kirwan's results in this latter direction: 
%namely, we still work within the context of symplectic geometry, but now consider
%the $K$-theory of the symplectic manifolds in question. 

Our primary technical tool is the $K$-theoretic version of the
Atiyah-Bott lemma, which states that the equivariant Euler class of a
$G$-bundle is not a zero-divisor, provided that a circle subgroup
$S^{1}\subset G$ fixes the base but acts non-trivially away from the
zero section.  If we use the
norm square of the moment map to build our manifold out of equivariant
strata, then the normal bundles satisfy this Atiyah-Bott condition, and 
we can then show that the
norm square of the moment map is an equivariantly perfect Morse-Kirwan
function. In other words, we can compute the equivariant $K$-theory of
the whole space entirely in terms of the $K$-theory of the critical
sets (note that the periodicity of complex $K$-theory eliminates the
degree shifts that appear in the cohomology version).

A proof of the Atiyah-Bott lemma in the context of algebraic
$K$-theory, for abelian actions, is given in \cite{VV03}. We give a
topological proof based on that of \cite{VV03}, extending this
result to 
actions of compact connected Lie
groups, as well as to bundles possessing only a $\Spinc$-structure instead
of a complex structure. 
It turns out that the $K$-theoretic Atiyah-Bott lemma is
stronger than its counterpart in rational cohomology; we do not need
extra assumptions on the torsion of the space or of the group. We take
some care to discuss these differences, and indicate how we
expect this $K$-theoretic Atiyah-Bott lemma, which is central to many
results in the theory of Hamiltonian $G$-manifolds, to give new information in the
context of symplectic geometry.

Kirwan originally proved surjectivity in the context of rational
cohomology. Our $K$-theoretic adaptation of Kirwan's work has a
slightly different flavor, with subtle differences arising from the
torsion components. For example, we are no longer simply computing
Betti numbers or Poincar\'e polynomials, and group extensions are not
necessarily simply direct sums. If we were to eliminate the torsion in $K$-theory
by tensoring with $\Q$, we recall from \cite{Ati-Hir} that the Chern
character gives an isomorphism between rational $K$-theory and
rational cohomology. In this sense, our Theorem~\ref{thm:KthyKirwan} can be viewed as an
extension of Kirwan's results for rational cohomology. Indeed, we
believe this $K$-theoretic extension to be more natural than a
corresponding extension to integral cohomology. As we mentioned above, 
when working with
integral cohomology, one needs to place additional constraints on the
torsion to establish the Atiyah-Bott lemma and its consequences (see
\cite{TW03}), but such torsion constraints are not present in the
$K$-theoretic version.

In Section \ref{formality}, we discuss equivariant formality. While the usual notion of equivariant formality does not have a counterpart in equivariant $K$-theory, we define a slightly weaker version of equivariant formality and prove in  Theorem~\ref{theorem:formality} that it holds for all Hamiltonian $G$-spaces. Our primary tool here is the equivariant K\"unneth theorem spectral sequence of Hodgkin, as refined by Snaith and McLeod (see \cite{Hod75}, \cite{Sna72}, and \cite{McL79} respectively).
As a consequence of equivariant formality we obtain a different sort of surjectivity, in this case for the forgetful homomorphism
\begin{equation}\label{eq:forgetful-k}
K_{G}^{*}(M) \twoheadrightarrow K^{*}(M)
\end{equation}
taking equivariant $K$-theory to non-equivariant $K$-theory. At the level of bundles, this implies that every complex vector bundle over a Hamiltonian $G$-space admits a stable equivariant structure. We go on to show in Theorem~\ref{equivariant-line-bundles}
that we can remove the word \emph{stable} for the case of complex line bundles,
%using cohomological techniques to
establishing the surjectivity of the forgetful map
\begin{equation}\label{eq:forgetful-pic}
\Pic_{G}(M) \twoheadrightarrow \Pic(M)
\end{equation}
for the Picard group of isomorphism classes of (equivariant) complex line bundles.

We translate our statements about line bundles and the Picard group into statements about integral cohomology in Appendix~\ref{appendix-line-bundle}. In particular, we show that the forgetful map
\begin{equation}\label{eq:forgetful-h2}
H_{G}^{2}(M;\Z) \twoheadrightarrow H^{2}(M;\Z)
\end{equation}
is surjective. For this result,
% which in the $K$-theory setting is the statement that every
%vector bundle can be given an equivariant structure. 
%We conjecture
%this to be true for compact Hamiltonian $G$-spaces and prove a partial
%result for line bundles in the case where the fundamental group of $G$
%has no torsion. 
%In the process,
we use the classification of
$G$-equivariant line bundles in terms of the equivariant cohomology
$H^{2}_{G}(M;\Z)$. This fact is non-trivial, and although it appears
in various places throughout the literature, we offer a previously
unpublished elementary proof due to Peter Landweber.
% in Appendix~\ref{appendix-line-bundle}.
At the core of this argument is
the fact that the homomorphism $\alpha : K_{G}(M) \to K(M_{G})$
retracts onto a map $\Pic_{G}(M) \to \Pic(M_{G})$ of Picard groups. A
technical lemma requiring only linear algebra then shows that this
retraction extends to the completion of $K_{G}(M)$ at the augmentation
ideal, and thus the map of Picard groups is in fact a retract of
the Atiyah-Segal isomorphism $\widehat{\alpha} : K_{G}(M)^{\wedge} \to
K(M_{G})$ (see \cite{Ati-Seg}). Identifying $\Pic(M_{G})$ with
$H^{2}_{G}(M;\Z)$ gives the desired classification of $G$-line
bundles.
%As a further consequence, we observe that $K_{G}(M)$ retracts
%onto $H^{2}_{G}(M;\Z)$, which can be used to given an alternative proof of the surjectivity of (\ref{eq:forgetful-h2}), derived from the surjectivity of (\ref{eq:forgetful-k}).

%motivating our use of cohomological techniques
%to study equivariant $K$-theory.
%At the core of this argument is the fact that the homomorphism $\Pic_{G}(M) \to \Pic(M_{G})$ of Picard groups is a retract of the Atiyah-Segal homomorphism $\alpha : K_{G}(M) \to K(M_{G})$ from equivariant $K$-theory the $K$-theory of the Borel construction $M_{G} = EG \times_{G} M$. Taking the completion of $K_{G}(M)$ with respect to the augmentation ideal, a linear algebra computation shows that the homomorphism on Picard groups is in fact the retraction of the Atiyah-Segal isomorphism $\widehat{\alpha} : K_{G}(M)^{\wedge} \to K(M_{G})$.
%, taking advantage of the intersection of $K$-theory and cohomology theory.
% GDL 3/11/05

We conclude with a discussion of Morse-Kirwan functions in
Appendix~\ref{morse-kirwan}, showing that the results of Morse-Bott
theory still hold in $K$-theory for such minimally degenerate
functions, provided that we work in terms of strata rather than simply
critical sets. Such an argument is not strictly necessary when working
with the norm square of the moment map, as it is better behaved than a
generic Morse-Kirwan function (see \cite{Ler04,Woo02}). However, we
offer it for the sake of completeness, and in case the reader is
interested in Morse-Kirwan functions other than the norm square of the
moment map.

The research presented in this paper represents an ongoing
project. 
Beyond the fundamental result of Kirwan surjectivity, 
another natural step 
is to perform computations of the kernel of the Kirwan map as in \cite{Gol02,TW03} 
to explicitly compute 
the $K$-theory of symplectic quotients. This also involves 
expected $K$-theoretic analogues of results of
Martin \cite{Mar00}, in order to 
relate the $K$-theory of the symplectic quotient by a nonabelian
Lie group $G$ to that of the corresponding 
abelian symplectic quotient by a maximal torus $T \subseteq
G$. Moreover, 
now that we have extended Kirwan's results to complex $K$-theory, it
is natural to ask whether there are analogues for other generalized
cohomology theories, such as real $K$-theory, cobordism, or elliptic
cohomology. We believe that the Atiyah-Bott lemma holds for
equivariant cobordism, as described in \cite{Sin01}.
One can also ask whether there are $K$-theoretic extensions of
similar cohomological results for other geometries possessing moment
maps, such as hyperk\"ahler, contact, Sasakian, or 3-Sasakian
geometry.

The authors would like to thank Matthias Franz, Jean-Claude Hausmann,
Tara Holm, Mike Hopkins, Lisa Jeffrey, Allen Knutson, Peter Landweber,
Haynes Miller, and Jonathan Weitsman for their insight and many
helpful discussions. 
We would also like to thank the American
Institute of Mathematics for hosting a workshop on the subject of Kirwan 
surjectivity, which originally inspired this project, as well as the Banff International Research Station, at which the authors presented these results to the symplectic geometry community. The second author would like to thank the Fields Institute and the University of Toronto for their support and hospitality while pursuing
this research and preparing this manuscript.

\section{The $K$-theoretic Atiyah-Bott lemma}\label{sec:AtiyahBott}

The main tool used in our proof of the surjectivity
theorem is Lemma~\ref{lemma:AB} below, which is the
$K$-theoretic analogue of the key fact behind many results in
equivariant symplectic geometry. The original version of the lemma proved by Atiyah and Bott is \cite[Proposition 13.4, p.~606]{AB82}, and is stated in 
terms of equivariant cohomology. An algebraic $K$-theory version of this lemma, for torus actions (or more precisely, for actions of diagonalizable group schemes of finite type on separated noetherian regular algebraic spaces) can be found in \cite[\S 4]{VV03}. Our proof is modeled on this algebraic $K$-theory proof, but we
work in the topological context, and we extend the result to general (non-abelian) compact connected Lie group actions on $\Spinc$-bundles.
First, we require a simple technical lemma, which allows us to impose complex structures on the bundles involved.

\begin{lemma}\label{lemma:complex-structure}
Let $G = H \times S^{1}$ be a compact Lie group, and let $E \to X$ be a real $G$-vector bundle over a $G$-space $X$. Suppose that the action of the $S^{1}$-factor of $G$ fixes precisely the zero section $X \subset E$. Then $E$ admits a $G$-invariant complex structure (and is therefore orientable).
\end{lemma}

\begin{proof}
%If $E$ and $X$ are smooth, then we can identify $E$ with the normal bundle for the $H$-equivariant inclusion of the $S^{1}$-fixed point set $X = E^{S^{1}} \hookrightarrow E$ via the zero section. Such a normal bundle necessarily admits an $S^{1}$-invariant complex structure, as we show below. Since the $H$-action commutes with the $S^{1}$-action, the complex structure on $E$ is $H$-invariant, and thus $G = H \times S^{1}$-invariant.
A real vector bundle $E$ can be viewed as a complex vector bundle $E\otimes_{\R} \C$ together with an involution corresponding to conjugation on the $\C$-factor. We recall that a complex $H\times S^{1}$-bundle over a space with trivial $S^{1}$-action can be decomposed into isotypic pieces, yielding an isomorphism of complex $G$-bundles
\begin{equation*}
E \otimes_{\R} \C 
\cong {\bigoplus}_{\lambda\in\Z} E_{\lambda} \otimes_{\C} \C_{\lambda},
\end{equation*}
where the $E_{\lambda}$ are complex $H$-bundles, and $\C_{\lambda}$ is the irreducible representation of $S^{1}$ on $\C$ with $z\in S^{1}$ acting as multiplication by $z^{\lambda}$. In our case, there is no trivial component with $\lambda = 0$ since $S^{1}$ fixes only the zero section. Furthermore, since $E \otimes_{\R} \C$ is self-conjugate, we have $E_{-\lambda} \cong \overline{E}_{\lambda}$ as complex bundles, and we will show that we can identify $E$ with the complex bundle $\bigoplus_{\lambda > 0} E_{\lambda} \otimes_{\C} \C_{\lambda}$ as real $G$-bundles.

The non-trivial irreducible real representations of $S^{1}$ are of the form $\R^{2}_{\lambda}$ for integers $\lambda > 0$, where
$\R^{2}_{\lambda} \otimes_{\R} \C \cong \C_{\lambda} \oplus \C_{-\lambda}$. Such a representation is called \emph{complex}, as the commuting-field of complex $S^{1}$-equivariant endomorphisms of $\R^{2}_{\lambda} \otimes_{\R} \C$ is $2$-dimensional and can be identified with $\C \otimes_{\R} \C \cong \C \oplus \overline{\C}$.  Continuing to work with the complexifications, we have the decomposition into isotypic components described in \cite[\S 8]{Ati-Seg} and \cite[\S 2]{Seg68},
$$E \otimes_{\R} \C \cong {\bigoplus}_{\lambda > 0}
\Hom_{\C}^{S^{1}} \left( \R^{2}_{\lambda} \otimes_{\R} \C, \, E \otimes_{\R}\C \right) \otimes_{\End_{\C}^{S^{1}}(\R^{2}_{\lambda}\otimes_{\R}\C)} \bigl( \R^{2}_{\lambda} \otimes_{\R} \C\bigr),$$
which respects the conjugation involution. Restricting to the real parts fixed under conjugation, we have the corresponding decomposition
$$ E \cong {\bigoplus}_{\lambda > 0}\Hom_{\R}^{S^{1}}\left(\R^{2}_{\lambda}, E\right) \otimes_{\End_{\R}^{S^{1}}(\R^{2}_{\lambda})} \R^{2}_{\lambda}.$$
Finally, identifying $\R^{2}_{\lambda} \cong \C_{\lambda}$ for all integers
$\lambda > 0$ gives us the identification $\End_{\R}^{S^{1}}(\C_{\lambda}) \cong \C$ for the commuting-field. We can therefore make the identification
\begin{equation}\label{eq:complex-structure}
E \cong {\bigoplus}_{\lambda > 0} \Hom_{\R}^{S^{1}}\left(\C_{\lambda}, E\right)\otimes_{\C} \C_{\lambda},
\end{equation}
where the $S^{1}$-invariant complex structure on each component comes from the $\C_{\lambda}$-factors. The action of the group $H$ on each component comes from the first factor, with $H$ acting trivially on the $\C_{\lambda}$-factors. The complex structure is therefore $H$-invariant as well.

If $E$ and $X$ are smooth, we can construct this complex structure explicitly
from the vector field $\xi$ on $E$ generating the $S^{1}$-action. Since the $S^{1}$-action maps each fiber onto itself, the infinitesimal action $\xi$ is a vertical vector field. Furthermore, since $S^{1}$ acts linearly on the fibers, the vector field $\xi$ is given by fiberwise multiplication by an endomorphism $I_{x} : E_{x} \to E_{x}$ for each point $x\in X$ in the base. 
On the isotypic component of $E$ transforming like $\R^{2}_{\lambda}$, we can then choose our complex structure to be the rescaling $J = \frac{1}{\lambda} I$, which yields the identification $\R^{2}_{\lambda}\cong \C_{\lambda}$ used above.
\end{proof}

\begin{remark}\label{remark:normal}
If in addition $E \to X$ is a smooth vector bundle, then the conditions of Lemma~\ref{lemma:complex-structure} allow us to identify $E$ with the normal bundle for the $H$-equivariant inclusion of the $S^{1}$ fixed point set $X = E^{S^{1}}\hookrightarrow E$ via the zero section. It is well known that the normal bundle to a connected component of the fixed point set for an $S^{1}$-action admits a complex structure, essentially by the proof given above, and Lemma~\ref{lemma:complex-structure} simply generalizes this fact.
%In our case, the complex structure is also $H$-invariant since the $H$-action commutes with the $S^{1}$-action.
\end{remark}

We are now ready to prove our $K$-theory analogue of the Atiyah-Bott lemma.

\begin{lemma}\label{lemma:AB}
Let a compact connected Lie group $G$ act fiberwise linearly on a $\Spinc$-vector bundle \(\pi: E \to X\) over a compact $G$-manifold $X$. 
Assume that a circle subgroup $S^1\subset G$ acts by restriction on
$E$ so that its fixed point set is precisely the zero section $X$.
Choose an invariant metric
on $E$ and let $D$ and $S$ denote the disc and sphere bundles,
respectively. Then the long exact sequence for the pair $(D,S)$
in equivariant $K$-theory splits into short exact sequences
\[
\xymatrix{
0 \ar[r] & K^*_{G}(D, S) \ar[r] & K^*_{G}(D) \ar[r] & K^*_{G}(S)
\ar[r] & 0.
}
\]
\end{lemma}

\begin{proof}
We will show that the first homomorphism $K^{*}_{G}(D,S)\to K^{*}_{G}(D)$ is injective. By the equivariant Thom isomorphism (see \cite{Ati68}),
%\cite{Kar70}, \cite{Seg68}),
we have for the domain $K^{*}_{G}(D,S) \cong K^{*}_{G}(X)$,
and since the disc bundle $D$ retracts
equivariantly onto the zero section $X$,
we likewise have for the codomain
$K^{*}_{G}(D) \cong K^{*}_{G}(X)$. In terms of these
isomorphisms, as shown in the following commutative diagram
$$\xymatrix{
 	\cdots \ar[r] &
	K_{G}^{*}(D,S) \ar[r]  &
	K_{G}^{*}(D) \ar[r] \ar[d]^{\cong} &
	K_{G}^{*}(S) \ar[r] & \cdots \\
	& K_{G}^{*}(X) \ar[u]^{\cong}_{\mathrm{Thom}} \ar[r]^{\otimes\, e_{G}(E)}
	& K_{G}^{*}(X)
}$$
this first map is multiplication by the $G$-equivariant Euler class,
\begin{equation*}
e_{G}(E) = [\bbS^{+}_{E}] - [\bbS^{-}_{E}] \in K_{G}(X),
%\lambda_{-1}(E^{*}) = {\sum}_{p} (-1)^{p}\bigl[\Lambda_{\C}^{p}(E^{*})\bigr]\in K_{G}(X).
\end{equation*}
where $\bbS^{+}_{E}$  and $\bbS^{-}_{E}$ are the two complex half-spin representations determined by the $\Spinc$-structure on $E$ (note that it follows from the conditions of the lemma that $\rank E$ is even). If $E$ is in fact a complex bundle, then it admits a canonical $\Spinc$-structure with
half-spin representations
$$
\bbS^{+}_{E} = \Lambda^{\mathrm{even}}_{\C}(E^{*}), \qquad
\bbS^{-}_{E} = \Lambda^{\mathrm{odd}}_{\C}(E^{*}),
$$
and corresponding Euler class
\begin{equation}\label{eq:euler-class}
e_{G}(E) = \lambda_{-1}(E^{*}) = {\sum}_{p} (-1)^{p}\left[\Lambda_{\C}^{p}(E^{*})\right].
\end{equation}
See \cite{Bott69} for a discussion of the $K$-theory Euler class in the non-equivariant case.
To establish the injectivity of this map, we will demonstrate that this $K$-theoretic Euler class is not a zero-divisor. (We note that a related result concerning the invertibility of the equivariant $K$-theory Euler class appears as \cite[Lemma 2.7]{AS68II}, dating back to the early days of $K$-theory.)

Let $T$ be a maximal torus of $G$ containing the $S^{1}$-subgroup fixing $X$. Given any $G$-bundle over $X$, by restricting the action from $G$ to $T$ we obtain a $T$-bundle, giving us the forgetful map
$K_{G}^{*}(X) \to K_{T}^{*}(X)$. This homomorphism is injective by \cite[Proposition 4.9]{Ati68}, in which Atiyah argues that
% this map has a left inverse, a 
 the $K$-theory pushforward map defined in terms of the Dolbeault complex for $G/T$ gives a left inverse. The image of the forgetful map is precisely the Weyl invariants in $K_{T}^{*}(X)$. By the naturality of the Euler class, the image of $e_{G}(E)\in K_{G}(X)$ is $e_{T}(E)\in K_{T}(X)$, and so we can 
therefore restrict our attention to $T$ and show that $e_{T}(E)$ is not
a zero-divisor.

Factoring out our designated circle, the torus decomposes as $T \cong T' \times S^{1}$, where $T' = T/S^{1}$. Since the $S^{1}$-action fixes $X$, we obtain a $\Z_{2}$-graded ring isomorphism
\begin{equation}\label{eq:KT-factor}
	K_{T}^{*}(X) \cong K_{T'}^{*}(X) \otimes K_{S^{1}}^{*}(\pt)
	\cong K_{T'}^{*}(X)[ z, z^{-1} ]
\end{equation}
(see \cite{Seg68}), where $K^{0}_{S^{1}}(\mathrm{pt}) \cong R(S^{1}) \cong \Z[z,z^{-1}]$ is the representation ring of $S^1$, and
$K^{1}_{S^{1}}(\mathrm{pt}) = 0$ vanishes. Indeed, any
complex $S^{1}$-bundle $V$ over $X$ decomposes in terms of the $S^{1}$-action
as a finite sum $V \cong \bigoplus_{\lambda\in\Z}V_{\lambda}\otimes
\C_{\lambda}$, where the $V_{\lambda}$ are complex $T'$-vector bundles over $X$, and $\C_{\lambda}$ is the representation of $S^{1}$ on $\C$ with $z\in
S^{1}$ acting as multiplication by $z^{\lambda}$.

By Lemma~\ref{lemma:complex-structure}, the bundle $E$ admits a $T$-invariant complex structure, which by \eqref{eq:complex-structure} decomposes into isotypic components for the $S^{1}$-action as $E \cong \bigoplus_{\lambda > 0} E_{\lambda}\otimes \C_{\lambda}$.
% since $S^{1}$ fixes only the zero section $X$.
This complex structure induces a canonical $\Spinc$-structure on $E$, whose corresponding half-spin representations, and consequently its Euler class (\ref{eq:euler-class}), differ from that of the given $\Spinc$ structure (or the given complex structure) by a line bundle.
The $K$-theory class of this line bundle is of the form $[L]z^{k} \in K^{*}_{T'}(X)[z,z^{-1}]$, for some complex $T'$-equivariant line bundle $L$ and integer $k$.  Recalling that the Euler class (\ref{eq:euler-class}) for complex bundles is multiplicative, we compute
\begin{align}
\notag
e_{T}(E) &= [L]z^{k} \, e_{T} \left( {\bigoplus}_{\lambda > 0}E_{\lambda}\otimes \C_{\lambda} \right)
= [L]z^{k} \,{\prod}_{\lambda > 0} e_{T}(E_{\lambda} \otimes \C_{\lambda}) \\
\label{eq:euler-factors}
&= [L]z^{k}\,{\prod}_{\lambda > 0} \left( 1 - E_{i}^{*} z^{-\lambda} + \cdots \pm \Lambda^{\rank E_{\lambda}}(E_{\lambda}^{*})\, z^{-\lambda\rank E_{\lambda}} \right).
\end{align}
%By our assumption that the $S^{1}$-action leaves only $0\in E_{\mathrm{pt}}$ fixed, we have $\lambda_{i}\neq 0$ for all $i$, and thus $i^{*}e_{S^{1}}(E)\neq 0$.
%where the factor of $[L]z^{k}$, for some complex $T'$-equivariant line bundle $L$ and integer $k$, measures the difference between the given $\Spinc$-structure on $E$ and the $\Spinc$-structure corresponding to the complex structure.
To show that $e_{T}(E)$ is not a zero-divisor in $K_{T}^{*}(X)$, we show that each of its factors is not a zero-divisor. First, we recall that the class $[L] z^{k}\in K_{T}^{*}(X)$ of a complex line bundle is a unit in $K$-theory and thus not a zero-divisor. Secondly, a polynomial $p(z) = 1 + a_{1} z^{-\lambda} + \cdots + a_{n} z^{-n\lambda}$ for $\lambda > 0$ with constant term $1$ and coefficients $a_{i} \in K^{*}_{T'}(X)$ cannot be a zero-divisor in $K_{T}^{*}(X)$. Explicitly, suppose that $p(z) q(z) = 0$ for a finite Laurent series $q(z) = \sum_{i}b_{i}z^{i} \in K^{*}_{T'}(X)[z,z^{-1}]$. Then,
\begin{align*}
0
&= {\sum}_{i} \left( b_{i} z^{i} + a_{1}b_{i} z^{i - \lambda} + \cdots + a_{n} b_{i}z^{i-n\lambda} \right) \\
&= {\sum}_{i}\left( b_{i} + a_{1}b_{i+\lambda} + \cdots + a_{n}b_{i+n\lambda} \right) z^{i},
\end{align*}
and thus we can write each coefficient of $q(z)$ as
$$b_{i} = -a_{1}b_{i+\lambda} - \cdots - a_{n}b_{i+n\lambda}$$
in terms of higher degree coefficients.
% depending on whether $\lambda > 0$ or $\lambda < 0$ respectively. 
However, if $q(z) = \sum_{i}b_{i}z^{i}$ is a finite Laurent series, this immediately implies that $b_{i} = 0$ for all $i$, and thus $q(z) = 0$.
\end{proof}

\begin{remark}
The proof of the equivariant cohomology version of this lemma given by Atiyah and Bott in \cite{AB82} uses the filtration $F_{p} = H^{\geq p}_{T'}(X) \otimes
H^{*}_{S^{1}}(\pt)$ of $H^{*}_{T}(X)$ by the cohomological degree of the $H^{*}_{T'}(X)$-factor. With respect to this filtration, the projection onto the degree 0 component of the associated graded algebra is the restriction map
$H^{*}_{T}(X) \to H^{*}_{S^{1}}(\pt)$.
Acting on the associated graded algebra, multiplication by the equivariant Euler class descends to multiplication by its image in $H^{*}_{S^{1}}(\pt)$. By this argument, all that is necessary is to show that the equivariant Euler class is not a zero-divisor over a point.

We can reprise this proof in $K$-theory using the filtration of $K^{*}_{T}(X)$ by powers of the kernel of the restriction map $K^{*}_{T}(X) \to K^{*}_{S^{1}}(\pt).$
For $T = S^{1}$, this is the filtration
$F_{p} =[ \tilde{K}^{*}(X) ]^{p} \otimes K^{*}_{S^{1}}(\pt)$
using powers of the reduced $K$-theory in place of the cohomological degree.
%(Alternatively, we could use the Atiyah-Hirzebruch filtration on $K^{*}(X)$, or if $K^{*}(X)$ is torsion-free, then we could pull back via the Chern character the filtration on $H^{*}(X;\Q)$ by the cohomological degree.)
Once again, the problem is reduced to showing that the equivariant Euler class is not a zero-divisor over a point. A single fiber $E_{\pt}$ is a finite dimensional representation of $S^{1}$, which splits as a direct sum $E_{\pt} \cong \bigoplus_{i} \C_{\lambda_{i}}$ of one-dimensional irreducible representations, and \eqref{eq:euler-factors} becomes
\begin{equation*}%\label{eq:alpha}
%i^{*}e_{S^{1}}(E)  = 
e_{S^{1}}(E_{\mathrm{pt}}) 
 = z^{k}\,e_{S^{1}}\Bigl( {\bigoplus}_{i} \C_{\lambda_{i}} \Bigr)
 = z^{k}\,{\prod}_{i} e_{S^{1}}\bigl( \C_{\lambda_{i}} \bigr)
 = z^{k}\,{\prod}_{i} \bigl( 1 - z^{-\lambda_{i}} \bigr).
\end{equation*}
%where $z^{k}$ measures the difference, over a point, between the given $\Spinc$-structure and the canonical $\Spinc$-structure induced by the complex structure.
Since the Euler class over a point is a finite Laurent polynomial whose lowest (and highest) degree term has coefficient $\pm 1$, it is not a zero-divisor in $\Z[z,z^{-1}]$. Furthermore, this Euler class is not divisible by any prime $p\in \Z$ and is thus a primitive element of $K^{*}_{S^{1}}(\pt)$.
\end{remark}

The original Atiyah-Bott lemma in \cite{AB82} for cohomology is actually weaker than our result for $K$-theory. In Lemma~\ref{lemma:AB}, we
need not make any assumptions on the torsion of the $K$-theory
$K^{*}(X)$, since the $K$-theoretic Euler class restricted to a point is always primitive. 
In contrast, the cohomology Euler class may not
be primitive, because $e_{S^{1}}(E_{\mathrm{pt}}) = \prod_{i} \lambda_{i}u
\in H^{*}_{S^{1}}$, where $u$ is the generator of $H^{2}_{S^{1}}$, may
be divisible by primes $p\in \Z$. Therefore, to state the
analogous theorem for equivariant cohomology, we either require that
$H^{*}(X;\Z)$ be torsion-free, or we explicitly kill the torsion by
working with $H^{*}(X;\Q)$.

Atiyah and Bott also require a second condition on the torsion. In order to establish the injection $H^*_{G}(X;\Z) \to H^*_{T}(X;\Z)$, they require that the
cohomology $H^{*}(G;\Z)$ of the Lie group $G$ be torsion-free in order that 
$H^{*}(BG;\Z)$ be torsion-free and the fibration \(G/T \to BT \to BG\) behave like a product for
integral cohomology, i.e., $H^{*}(BT;\Z) \cong H^{*}(BG;\Z) \otimes H^{*}(G/T;\Z)$. However, %as the $K$-theory of a compact Lie
as the representation ring \(K_{G}^*\) is always torsion-free, %(see \cite{Ati65}),
there is no need for this condition when working with $K$-theory. Furthermore,
by Hodgkin's Theorem (see \cite{Hod67,McL79}), if a compact connected Lie group has no torsion in its fundamental group, then its $K$-theory is automatically torsion free.

In light of these observations, our $K$-theoretic Atiyah-Bott lemma is
actually stronger than its cohomological counterpart. It works without
localization or the need to tensor with the rationals, and it also works
for all compact manifolds $X$ and compact connected Lie groups $G$, not just
those whose cohomology is torsion-free. This is because Thom and Euler classes are better behaved in $K$-theory, allowing us to formulate our arguments with no mention of torsion. From another point of view,
recall that the Atiyah-Hirzebruch spectral sequence (see \cite{Ati-Hir}) implies that the order of the torsion subgroup in $K$-theory is at most equal to the order of the torsion subgroup in integral cohomology. We then see that passing from integral cohomology to $K$-theory eliminates just enough torsion for the Atiyah-Bott lemma to work.

%\smallskip

%Finally, we observe that Lemma~\ref{lemma:AB} holds under the slightly weaker
%condition that the bundle $E$ possess a $G$-invariant
%$\Spinc$-structure, as opposed to a complex structure. With such a
%structure, the Thom isomorphism still holds in equivariant $K$-theory,
%and the Euler class is given by the direct difference of the two
%half-spin bundles associated to $E$.

\section{Kirwan surjectivity in $K$-theory}\label{sec:surjectivity}

%\todo{Do we want $G$ to be connected? I think not}

The goal of this section is to prove a $K$-theoretic extension of the
surjectivity theorem of Kirwan. We now briefly recall the setting of
these results. Let $G$ be a
compact connected Lie group.
A \emph{Hamiltonian $G$-space} is a symplectic manifold $(M,\omega)$
on which $G$ acts by symplectomorphisms, together with a $G$-equivariant
moment map $\mu: M \to \g^*$ satisfying
Hamilton's equation
\[
\left< d\mu, X \right> = \imath_{X^{\sharp}}\omega,  \quad \forall X
\in \g,
\]
where \(X^{\sharp}\) denotes the vector field on $M$ generated by $X
\in \g$. We will assume throughout that the moment map $\mu: M \to
\g^*$ is proper, i.e., the preimage of a compact set is compact. 
We further assume that $0$ is a regular value so that
$G$ acts locally freely on the
level set $\mu^{-1}(0)$, i.e., the $G$-action has finite stabilizers. In this situation the {\em symplectic
quotient} or Marsden-Weinstein reduction of $M$ at $0$ is
$$ M \mod G := \mu^{-1}(0) \,/\, G,$$
viewed as the standard quotient if the action is free, or an orbifold otherwise.

The surjectivity theorem of Kirwan
gives a method of computing the rational cohomology ring of this
symplectic quotient using that of the original symplectic manifold
$M$ (see \cite{Kir84}). We do the same here, except in $K$-theory. Our main theorem is as
follows. 

\begin{theorem}\label{thm:KthyKirwan}
Let \((M,\omega)\) be a Hamiltonian $G$-space with proper moment map
$\mu: M \to \g^*$. 
Assume that $0$ is a regular value of $\mu$ so that the group $G$ acts locally freely on $\mu^{-1}(0)$. Then the Kirwan map $\kappa$
induced by the inclusion \(\iota : \mu^{-1}(0) \into M,\)
\[
\xymatrix{
K^*_G(M) \ar[r]^-{\iota^*} \ar[rd]_{\kappa} & K^*_G\bigl(\mu^{-1}(0)\bigr) \ar[d]^{\cong}\\
& K^*(M \mod G)
}
\]
is a surjection.
\end{theorem}

If the $G$-action on the level set $\mu^{-1}(0)$ is not free, then we take $K^{*}(M\mod G)$ to be the orbifold $K$-theory, which for a global quotient is given, and in some cases defined, by
\begin{equation}\label{eq:orbifold}
K_{\text{orb}}^{*}\bigl(\mu^{-1}(0) \,/\, G\bigr) \cong K_{G}^{*}\bigl(\mu^{-1}(0)\bigr).
\end{equation}
See, for example, \cite{AR03} for further discussion of orbifold $K$-theory.

Our strategy for the proof of Theorem~\ref{thm:KthyKirwan} is to
follow Kirwan's original proof of surjectivity for rational cohomology
\cite{Kir84} and verify that the cohomological aspects of her argument
extend to $K$-theory.  With this in mind, we now summarize Kirwan's
proof.  In \cite[\S 4]{Kir84}, Kirwan considers the norm square of the
moment map, whose minimum $\mu^{-1}(0)$ is automatically the level set
of 
interest, and shows that it is a minimally degenerate Morse function
(i.e., ``is Morse in the sense of Kirwan''). In \cite[\S 10]{Kir84},
she establishes that the results of Morse-Bott theory still hold in
this generalized setting.
%In Appendix~\ref{morse-kirwan}, we established the $K$-theoretic versions of the lemmas required to do this.
In \cite[\S 5]{Kir84}, Kirwan proves that this Morse function is
equivariantly perfect. She then uses this result to establish
surjectivity via an inductive argument using a stratification built
out of the Morse function. 

We begin with a few definitions. Let $M =  \bigsqcup_{\beta\in\mathcal{B}}S_{\beta}$ be a stratification
of $M$ indexed by a partially ordered set $\mathcal{B}$. Extending this
partial ordering to a total ordering, consider the unions
$$M_{<\beta} := {\bigsqcup}_{\gamma<\beta}S_{\gamma}, \qquad
  M_{\leq\beta} := {\bigsqcup}_{\gamma\leq\beta}S_{\gamma}, \qquad
  M_{\geq\beta} := {\bigsqcup}_{\gamma\geq\beta}S_{\gamma}.$$
The ordered collection of subsets $M_{\leq\beta}$ then gives a filtration
of $M$. In analogy to the Atiyah-Hirzebruch filtration (see \cite{Ati-Hir}), we obtain a filtration of $K^{*}(M)$ and a
spectral sequence satisfying
\begin{equation}\label{eq:filtration-sequence}
E_{1} = {\bigoplus}_{\beta\in\mathcal{B}} K^{*}\bigl(M_{\leq\beta},M_{<\beta}\bigr), \qquad
  E_{\infty} = \Gr\,K^{*}(M),
\end{equation}
which converges to the associated graded algebra of the $K$-theory of $M$.
Suppose further that this stratification is smooth, i.e., the strata $S_{\beta}$ are locally closed submanifolds of $M$, and that it satisfies the closure property 
\begin{equation}\label{eq:closureRelations}
\overline{S}_{\beta} \subseteq M_{\geq\beta},
\end{equation}
i.e., the closure of $S_{\beta}$ is contained in the strata above $S_{\beta}$ for each \(\beta \in {\cal B}\). If the stratification is smooth, then we have a normal bundle $N_{\beta}$ to $S_{\beta}$ in $M$, and due to the closure property we can choose a tubular neighborhood
$U_{\beta}$ of $S_{\beta}$ in $M$ diffeomorphic to $N_{\beta}$ such that
\begin{equation}\label{eq:tubular-neighborhood}
	U_{\beta} \subset M_{\leq\beta}, \qquad
 	U_{\beta}\setminus S_{\beta} \subset M_{<\beta}.
\end{equation}
Indeed, due to the closure property (\ref{eq:closureRelations}), any point sufficiently close to $S_{\beta}$ in a higher stratum must be in the closure of $S_{\beta}$, and hence must lie along a direction tangent to $S_{\beta}$.
The normal directions to $S_{\beta}$ must therefore initially extend into the lower strata. By excision, we therefore have
\begin{equation}\label{eq:excision}
K^{*}(M_{\leq\beta}, M_{<\beta}) \cong K^{*}(N_{\beta},N_{\beta}\setminus S_{\beta}).
% \cong K^{*}_{c}(N_{\beta}),
\end{equation}
(In terms of $K$-theory with compact supports, this is  $K^{*}_{c}(N_{\beta}) = \tilde{K}^{*}(N_{\beta}^{+})$.) If $N_{\beta}$ is complex (or $\Spinc$), then we have a Thom isomorphism
\begin{equation}\label{eq:thom-isomorphism}
%K^{*}_{c}(N_{\beta})
K^{*}(N_{\beta},N_{\beta}\setminus S_{\beta})
\cong K^{*-d(\beta)}(S_{\beta}),
\end{equation}
where the degree $d(\beta)$ of the stratum is the rank of its normal bundle $N_{\beta}$.
Combining the isomorphisms (\ref{eq:excision}) and (\ref{eq:thom-isomorphism}),
our spectral sequence (\ref{eq:filtration-sequence}) becomes
\begin{equation}\label{eq:stratification-sequence}
	E_{1} \cong {\bigoplus}_{\beta\in\mathcal{B}}K^{*-d(\beta)}(S_{\beta}), \qquad
	E_{\infty} = \Gr\,K^{*}(M).
\end{equation}
The same discussion holds for equivariant $K$-theory for a $G$-space with a
$G$-invariant stratification.

\begin{definition}
A smooth stratification $M =  \bigsqcup_{\beta\in\mathcal{B}}S_{\beta}$ is
called \emph{(equivariantly) perfect} for $K$-theory if the spectral sequence (\ref{eq:stratification-sequence}) for (equivariant) $K$-theory collapses at the $E_{1}$ page, or equivalently if we have short exact sequences of the form
\begin{equation}\label{eq:short-exact-sequence-s}
\xymatrix{
0 \ar[r] & K^{*-d(\beta)}\bigl(S_{\beta}\bigr) \ar[r]  & 
K^*\bigl(M_{\leq \beta}\bigr) \ar[r] & 
K^*\bigl(M_{<\beta}\bigr) \ar[r] & 0}
\end{equation}
for each $\beta\in\mathcal{B}$ (similarly for equivariant $K$-theory).
\end{definition}

In other words, whenever we add a stratum, the equivariant $K$-theory
of the new space including that stratum is an extension of, and
therefore surjects onto, the equivariant $K$-theory of the union of
all the lower strata. Although our primary goal in this section is to
prove a surjectivity result, we will instead attack this long exact
sequence from the opposite direction, proving that the leftmost maps of these exact
sequences are injective. We will do this using our results from Section~\ref{sec:AtiyahBott}.
%(\ref{eq:short-exact-sequence}) is injective.

We can refine this definition slightly if the stratification is
obtained from a proper Morse function (or more generally from a Morse-Bott or
Morse-Kirwan function).  Let $f : M \to \R$ be a proper Morse function on a
Riemannian manifold $M$, and consider its negative gradient flow. Let
$\{ C_{\beta} \}$ denote the connected components of the critical set
of $f$, and for each $C_{\beta}$ define the stratum $S_{\beta}$ to be
the set of points of $M$ which flow down to $C_{\beta}$ via their
paths of steepest descent. (Note that this is a decomposition of
$M$ by stable manifolds, not the unstable manifolds as is usually the
case in Morse theory.) The critical sets are partially ordered by $\beta \leq \gamma$
if $f(C_{\beta}) \leq f(C_{\gamma})$. Extending this to a total order, this gives us a smooth
stratification of $M$ satisfying the closure property
(\ref{eq:closureRelations}). Furthermore, the inclusion
$C_{\beta}\hookrightarrow S_{\beta}$ of the critical set into the
stratum induces an isomorphism $K^*(S_{\beta}) \cong K^*(C_{\beta})$. (For
Morse-Bott functions, there is a deformation retraction from the
stratum onto the critical set. For a general Morse-Kirwan function, we
obtain a retraction at the level of (equivariant) $K$-theory, as we
describe in Appendix~\ref{morse-kirwan}.) The spectral sequence
(\ref{eq:stratification-sequence}) then becomes
\begin{equation}\label{eq:morse-sequence}
	E_{1} \cong {\bigoplus}_{\beta\in\mathcal{B}}K^{*-d(\beta)}(C_{\beta}), \qquad
	E_{\infty} = \Gr\,K^{*}(M).
\end{equation}
In the equivariant case, we require that both the Morse function and the Riemannian metric be $G$-invariant, giving us a $G$-invariant stratification.

\begin{definition}
A Morse (or Morse-Bott or Morse-Kirwan) function $f: M \to \R$ on a Riemannian manifold is
called \emph{(equivariantly) perfect} for $K$-theory if the spectral sequence (\ref{eq:morse-sequence}) for (equivariant) $K$-theory collapses at the $E_{1}$ page, or equivalently if we have short exact sequences of the form
\begin{equation}\label{eq:short-exact-sequence-c}
\xymatrix{
0 \ar[r] & K^{*-d(\beta)}\bigl(C_{\beta}\bigr) \ar[r]  & 
K^*\bigl(M_{\leq \beta}\bigr) \ar[r] & 
K^*\bigl(M_{<\beta}\bigr) \ar[r] & 0}
\end{equation}
for each $\beta\in\mathcal{B}$ (and similarly for equivariant $K$-theory).
\end{definition}

\begin{remark}
We note that our definition of (equivariant) perfection given above
differs slightly from the usual one. 
The usual definition of (equivariant) perfection is in terms of
(equivariant) Poincar\'{e} polynomials---in particular, that the (equivariant) Morse
inequalities are in fact equalities. However, that definition is
concerned with only the ranks of the free components, as given by the
Betti numbers in rational or real cohomology, and not the torsion that
may be present in $K$-theory. Since we are here concerned with
statements that include the torsion components, we opt for the
definition in terms of the short exact sequences. 
\end{remark}

Before getting to our main results for this section, we first recall an elementary lemma, proved in \cite[\S 2]{Seg68} and \cite[p.~144]{May96}:

\begin{lemma}\label{lemma:bundle}
	Let $H\subset G$ be a Lie subgroup and let $X$ be a compact $H$-space.
	Then $G \times_{H} X = (G\times X)/H$ is a bundle over the homogeneous space $G/H$ with fiber $X$, and $K_{G}^{*}( G \times_{H} X) \cong K_{H}^{*}(X)$.
\end{lemma}

In the special case where $X$ is a point, then this is simply the
isomorphism $K_{G}(G/H) \cong K_{H}(\mathrm{pt}) = R(H)$. We are now
ready to show that the norm square of the moment map is equivariantly
perfect.

\begin{theorem}\label{thm:equivariantly-perfect}
Let $(M,\omega)$ be a Hamiltonian $G$-space with proper moment map
$\mu: M \to \g^*$.  Given a $G$-invariant inner product on the
coadjoint representation $\g^{*}$, the norm square of the moment map
$f = \| \mu \|^{2}$ is an equivariantly perfect Morse-Kirwan function
for $K$-theory. In addition, the degrees $d(\beta)$ of the critical
sets are all even and may thus be dropped when working with complex
$K$-theory.
\end{theorem}

\begin{proof}
%We now reproduce Kirwan's argument in further detail.
%Let $\g$ denote the Lie algebra of $G$, and choose a $G$-invariant inner product on $\g$. Then let \(f = \|\mu\|^2\) be the norm square of the moment map, which we use to construct a smooth
%$G$-invariant stratification of $M$.
%Let $\{ C_{\beta} \}$, indexed by a finite set $\mathcal{B}$, be the critical sets of $f$.
Since the moment map $\mu : M \to \g^{*}$ is $G$-equivariant, the
moment map image of the critical set of $f := \| \mu \|^{2}$ must be
composed of complete coadjoint orbits. Since $G$ is compact, each
coadjoint orbit intersects the positive Weyl chamber (given a choice
of Cartan subalgebra and positive root system) exactly once, except
for those orbits which intersect the boundary of the Weyl chamber.  We
can therefore index the critical sets $\{C_{\beta}\}$ by a finite
subset $\mathcal{B} \subset \t_{+}$ of a fixed positive Weyl chamber,
defining $C_{\beta}$ for $\beta\in\mathcal{B}$ to be the set of all
critical points of $f$ with moment map image in the coadjoint orbit of
$\beta$. Then, for any $\beta\in\mathcal{B}$, the function $f$ takes
the constant value $f(C_{\beta}) = \| \beta \|^{2}$ on the critical
set $C_{\beta}$. We observe that the critical sets $C_{\beta}$ may not
be connected, as there may be several connected components of the
critical set of $f$ with moment map image in the same coadjoint
orbit.  (Actually, the different connected components of $C_{\beta}$
may have distinct degrees, so Kirwan further decomposes the critical
set as the disjoint union of $C_{\beta} = \bigsqcup_{d}C_{\beta,d}$
according to degree.) Since $\mu$ is proper, its norm square \(f:=\|\mu\|^2\) is
proper, and the critical sets $C_{\beta}$ are compact.

Choose a $G$-invariant metric on $M$ (which, together with the symplectic form, induces a compatible $G$-invariant almost complex structure on $M$).
%we can consider the negative gradient flow of the Morse function $f : M \to \R$.
%For each \(\beta \in {\cal B}\), we define the
%stratum $S_{\beta}$ to be the set of points of $M$ which flow down to the critical set
%$C_{\beta}$ via their paths of steepest descent. 
%The space $M$ then decomposes as the disjoint union \(M = \bigsqcup_{\beta \in 
%{\cal B}} S_{\beta}\) of these strata.
%(We note that this is a decomposition of $M$ by stable manifolds,
%not the unstable
%manifolds as is usually the case in Morse theory.)
%Furthermore, this decomposition is a stratification of $M$, in that the indexing set $\mathcal{B}$ has a strict partial order, here given by simply $\gamma > \beta$ if and only if $\| \gamma \|^{2} > \| \beta \|^{2}$, such that 
%\begin{equation}\label{eq:closureRelations}
%\overline{S}_{\beta} \subseteq {\bigsqcup}_{\gamma \geq \beta} S_{\gamma},
%\end{equation}
%for every \(\beta \in {\cal B}.\) Thus any stratum above $S_{\beta}$
%is contained in the closure of $S_{\beta}$. 
%In particular, the bottom
%stratum $S_{0}$ is open and dense in $M$ and contains the minimal critical set $C_{0} = \mu^{-1}(0)$.
Since
the norm square of the moment map is a $G$-invariant function,  this stratification is $G$-invariant. We can therefore consider the spectral
sequence~(\ref{eq:morse-sequence}) for equivariant $K$-theory.

Let $N_{\beta}$ be the normal bundle to $S_{\beta}$ in $M$, and let
$D(S_{\beta})$ and $S(S_{\beta})$ denote the disc and sphere bundles,
respectively, of $N_{\beta}$.
%By the tubular neighborhood theorem,
%there is a neighborhood
%$U_{\beta}$ of $S_{\beta}$ in $M$ diffeomorphic to $N_{\beta}$. Furthermore,
%we can choose this neighborhood $U_{\beta}$ such that
%any point in $U_{\beta} \setminus S_{\beta}$, i.e., any point
%in the normal bundle which is not in the zero section, is contained in
%some $S_{\gamma}$ for \(\gamma < \beta\). We can do this due to the
%closure relations in the stratification~\eqref{eq:closureRelations},
%as any point in a stratum above $\beta$ is in the closure of
%$S_{\beta}$ and hence must lie along a direction
%tangent to $S_{\beta}$. Thus the tubular neighborhood $U_{\beta}$ is
%contained entirely in \(\bigsqcup_{\gamma \leq \beta} S_{\gamma}.\)
We now claim
that there is a commutative diagram 
\begin{equation}\label{eq:commutative-diagram}
\xymatrix{
\cdots \ar[r] & K^*_G\bigl(M_{\leq \beta}, M_{<\beta}\bigr) \ar[r]^-{a} \ar[d]^{\cong} & 
K^*_G\bigl(M_{\leq \beta}\bigr) \ar[r] \ar[d]^{i^*} & 
K^*_G\bigl(M_{<\beta}\bigr) \ar[r] & \cdots \\
\cdots \ar[r] & K^*_G\bigl(D(S_{\beta}), S(S_{\beta})\bigr) \ar[r]^-{b} & K^*_G\bigl(D(S_{\beta})\bigr) \ar[r]
& K^*_G\bigl(S(S_{\beta})\bigr) \ar[r] & \cdots \\
}
\end{equation}
whose rows are long exact sequences.
The left vertical arrow is an isomorphism by excision (\ref{eq:excision}),
and the middle vertical arrow is the pullback
homomorphism induced by the inclusion \(i: D(S_{\beta}) \to
U_{\beta} \into M_{\leq\beta}\) of the tubular neighborhood $U_{\beta}$
from (\ref{eq:tubular-neighborhood}).
In order to prove that the top long exact sequence splits,
it is sufficient to show that the bottom long exact sequence splits.
Indeed, if the arrow $b$ injective, then 
the composition of the arrow $a$ with the pullback $i^{*}$ is injective,
and thus the arrow $a$ must be injective.

Recalling Lemma~\ref{lemma:AB}, in order to show that the long exact
sequence for the Thom complex $\bigl( D(S_{\beta}),S(S_{\beta})
\bigr)$ splits, we must find a circle subgroup of $G$ which acts on
the normal bundle $N_{\beta}$ fixing only the zero section. This does
not necessarily hold over the whole stratum $S_{\beta}$. However, the
negative gradient flow of the norm square of the moment map always
converges when $M$ is a Hamiltonian $G$-space with a proper moment map. 
(This result is due to an unpublished
manuscript of Duistermaat, whose argument is reproduced in
\cite{Ler04}. Another proof can be found in \cite[Appendix B]{Woo02}.)
The negative gradient flow therefore gives an equivariant deformation
retraction of the stratum $S_{\beta}$ onto the critical set
$C_{\beta}$. (Kirwan does not make use of this fact in \cite[\S
10]{Kir84}. Instead she establishes Morse theoretic results for minimally
degenerate functions which do not necessarily satisfy this special
property; we describe her argument, translated to the $K$-theory
setting, in
Appendix~\ref{morse-kirwan}.)  Furthermore, Kirwan observes in
\cite[\S 3]{Kir84} that each critical subset $C_{\beta}$ is of the
form \(G \times_{\Stab(\beta)} \bigl(Z_{\beta} \cap
\mu^{-1}(\beta)\bigr)\), where $Z_{\beta}$ is a critical set of the
non-degenerate Morse-Bott function $\mu^{\beta} = \langle \mu, \beta
\rangle$, the component of the moment map $\mu$ in the $\beta$
direction. Here $\Stab(\beta)\subset G$ is the isotopy subgroup
stabilizing the element $\beta\in \mathfrak{t}_{+}\subset\g$.  By
Lemma~\ref{lemma:bundle}, we have an isomorphism
\begin{equation}\label{eq:restriction}
K^*_G(C_{\beta}) \cong K^*_{\Stab(\beta)}\bigl(Z_{\beta} \cap \mu^{-1}(\beta)\bigr).
\end{equation}
Restricting the normal bundle $N_{\beta}$ first to the critical set
$C_{\beta}$ and then to the fiber $Z := Z_{\beta}\cap\mu^{-1}(\beta)$
over the identity coset of the homogeneous space $G/\Stab(\beta)$, let
$D(Z)$ and $S(Z)$ denote the disc and sphere bundles of
$N_{\beta}|_{Z}$.  The bottom long exact sequence of
(\ref{eq:commutative-diagram}) is therefore isomorphic to the long
exact sequence for the Thom complex $\bigl(D(Z),S(Z)\bigr)$ in
$\Stab(\beta)$-equivariant $K$-theory:
\begin{equation}\label{eq:LESforZ}
\xymatrix{
\cdots \ar[r] & K^*_{\Stab(\beta)}\bigl(D(Z), S(Z)\bigr) \ar[r] &
K^*_{\Stab(\beta)}\bigl(D(Z)\bigr) \ar[r] & K^*_{\Stab(\beta)}\bigl(S(Z)\bigr) \ar[r] & \cdots \\
}
\end{equation}

Let $T_{\beta}$ be the subtorus of $\Stab(\beta)$ generated by the element $\beta\in\g$. (If $\beta$ is an irrational element that does not generate a circle, then we let $T_{\beta}$ be the closure of the one parameter subgroup generated by $\beta$, which still gives us a torus.)
We recall that the critical set $Z_{\beta}$ of the moment map component $\mu^{\beta}$ is a (union of certain components of the) fixed point
set of $T_{\beta}$. It follows that this subtorus $T_{\beta} \subset \Stab(\beta)$ fixes $Z_{\beta}\cap \mu^{-1}(\beta)$ and acts on the normal bundle $N_{\beta}$ with no non-zero fixed vectors. 
Kirwan argues in \cite[\S 4]{Kir84} that the normal bundle $N_{\beta}$ is a complex vector bundle, with complex structure inherited from the $G$-invariant almost complex structure on $M$ induced by the symplectic form and metric. 
Therefore, by the $K$-theoretic Atiyah-Bott Lemma~\ref{lemma:AB}, the long exact sequence (\ref{eq:LESforZ}) splits into short exact sequences. It follows that the long exact sequences in (\ref{eq:commutative-diagram}) also split, which establishes the short exact sequences (\ref{eq:short-exact-sequence-c}). The norm square of the moment map is therefore equivariantly perfect in $K$-theory. Furthermore, as the normal bundle $N_{\beta}$ is complex, its real rank is even, and so the degrees $d(\beta)$ are even and do not affect complex $K$-theory.
\end{proof}
%\todo{Does this argument need to be refined, possibly using the isomorphism above?}

Now that we have established that the norm square of the moment map is equivariantly perfect, we can use the corresponding stratification to prove Kirwan surjectivity.

\medskip

\noindent\textbf{Proof of Theorem~\ref{thm:KthyKirwan}:} Finally, we
prove Kirwan surjectivity via an inductive argument.  Building the
space $M$ one stratum at a time, the short exact sequences
(\ref{eq:short-exact-sequence-c}) give us a chain of surjections
$$ \cdots \twoheadrightarrow K_{G}^{*}\bigl( M_{\leq\beta}\bigr)
\twoheadrightarrow K_{G}^{*}\bigl( M_{< \beta} \bigr)
\twoheadrightarrow
\cdots \twoheadrightarrow K_{G}^{*}(S_{0}),$$ which
ends with the bottom stratum $S_{0}$. The stratum $S_{0}$ is open and
dense in $M$ and retracts onto the minimal critical set $C_{0} =
\mu^{-1}(0)$.  Therefore \(K_{G}^{*}\bigl( M_{\leq \beta}\bigr)
\twoheadrightarrow K_{G}^{*}(S_0) \cong K_{G}^*(\mu^{-1}(0))\) for
every \(\beta \in {\cal B}.\) Since the maps \(K_{G}^{*} \bigl(
M_{\leq \beta}\bigr) \twoheadrightarrow K_{G}^{*}\bigl( M_{< \beta}
\bigr)\) are surjections for all $\beta$, there is no \(\lim^1\) term
and we may conclude that 
\[
K_{G}^{*}(M) = K_{G}^*\bigl(
\lim_{\longrightarrow}  M_{\leq \beta} \bigr) =
\lim_{\longleftarrow} 
K_{G}^{*}\bigl( M_{\leq \beta}\bigr).
\]
Hence the restriction map gives a surjection
$$\kappa: K_{G}^{*}(M) \twoheadrightarrow K_{G}^{*}\bigl(\mu^{-1}(0)\bigr).$$
Finally, if $0$ is a regular value, then $\mu^{-1}(0)$ is a manifold, and we
can compose this surjection with the isomorphism $K_{G}^{*}(\mu^{-1}(0)) \cong K^{*}( \mu^{-1}(0) / G)$, using orbifold $K$-theory as in (\ref{eq:orbifold}) if necessary.
\hfill$\Box$

\medskip

\section{Equivariant formality}\label{formality}

In this section we discuss the $K$-theoretic analogue of the Kirwan-Ginzburg equivariant formality theorem for compact Hamiltonian $G$-spaces. Recall that the equivariant cohomology of a $G$-space $M$ is defined in terms of the Borel construction by $H_{G}^{*}(M) := H^{*}(M_{G})$, where $M_{G} = M \times_{G} EG$. In \cite{GKM}, Goresky, Kottwitz, and MacPherson call a $G$-space $M$ \emph{equivariantly formal} if the Leray-Serre spectral sequence for the cohomology of the fibration $M \to M_{G} \to BG$ collapses at the $E_{2}$ page. When working with coefficients in a field, equivariant formality is equivalent to the existence of an isomorphism
\begin{equation}\label{eq:cohomology-formality}
H_{G}^{*}(M) \cong H^{*}(M) \otimes H_{G}^{*}(\pt)
\end{equation}
as modules over $H_{G}^{*}(\pt)$ or the corresponding statement about the Betti numbers or Poincar\'{e} polynomials.

Unfortunately, neither of these versions of equivariant formality applies to equivariant $K$-theory. Requiring an isomorphism analogous to (\ref{eq:cohomology-formality}) is too strong a condition due to the possible presence of torsion. On the other hand, a statement in terms of the Leray-Serre spectral sequence applies only to the Borel construction of equivariant $K$-theory, given by $\mathcal{K}_{G}^{*}(M) := K^{*}(M_{G})$. For the Atiyah-Segal equivariant bundle construction of $K_{G}^{*}(M)$, we do not have such a spectral sequence at our disposal. Instead, we consider the following weakened version of equivariant formality:

\begin{definition}
We call a $G$-space $M$ \emph{weakly equivariantly formal} for $K$-theory if the map
\begin{equation}\label{eq:weak-formality}
	K_{G}^{*}(M) \otimes_{R(G)} \Z \to K^{*}(M),
\end{equation}
induced by the forgetful homomorphism $K_{G}^{*}(M) \to K^{*}(M)$,
is an isomorphism.
\end{definition}

Here, the representation ring $R(G)$ acts on $\Z$ via the augmentation homomorphism $\epsilon : R(G) \to \Z$ taking a virtual representation to its dimension. So, rather than attempting to factor $R(G)$ out of the equivariant
$K$-theory, we instead divide by the $R(G)$-action, or more precisely by the action of the augmentation ideal $I(G) = \Ker \epsilon$.
The relationship between the forgetful homomorphism and the map (\ref{eq:weak-formality}) is given by the commutative diagram
$$\xymatrix{
	K_{G}^{*}(M) \ar@{->>}[d] \ar[dr] \\
	K_{G}^{*}(M) \otimes_{R(G)} \Z \ar[r] & K^{*}(M)
}$$
from which we obtain the following immediate consequences:

\begin{proposition}\label{prop:formality}
	If a $G$-space $M$ is weakly equivariantly formal for $K$-theory, 
	then the forgetful map $K_{G}^{*}(M) \to K^{*}(M)$ is surjective,
	and thus every complex vector bundle $E$ over $M$ admits a stable
	equivariant structure, i.e., the $G$-action on $M$ lifts to a
	fiberwise linear $G$-action on $E \oplus \underline{k}$ for some
	trivial bundle $\underline{k} = M \times  \C^{k}$. In addition,
	the kernel of the forgetful map is $I(G)\cdot K_{G}^{*}(M).$
\end{proposition}

We observe that our statement (\ref{eq:weak-formality}) of weak equivariant formality is closely related to the K\"unneth formula in equivariant $K$-theory. In non-equivariant $K$-theory, we have a K\"unneth formula (see \cite{Ati62}) which expresses the $K$-theory of a product $K^{*}(X\times Y)$ in terms of a short exact sequence involving the product of the $K$-theories $K^{*}(X)\otimes K^{*}(Y)$ together with torsion information. On the other hand, equivariant $K$-theory is a module over the representation ring $R(G)$, and for two $G$-spaces $X$ and $Y$, the 
equivariant $K$-theory of the product $K_{G}^{*}(X\times Y)$ can be expressed in terms of the tensor product $K_{G}^{*}(X)\otimes_{R(G)}K_{G}^{*}(Y)$
over $R(G)$, together with $R(G)$-torsion information.  In \cite{Hod75}, Hodgkin constructs a K\"unneth spectral sequence, related to an Eilenberg-Moore spectral sequence, starting with
$$E_{2}(X,Y) = \Tor^{*}_{R(G)}\bigl(K_{G}^{*}(X),K_{G}^{*}(Y)\bigr),$$
where the $\Tor$ groups are the derived functors of the tensor product
(see \cite{CE56} or any homological algebra textbook).
Here, we use Hodgkin's convention of labeling the higher $\Tor$ groups with a negative superscript in place of the more standard positive subscript, writing
$$\Tor^{-p}_{R(G)}(A,B) = \Tor_{p}^{R(G)}(A,B) \text{ with } p \geq 0,$$
which places Hodgkin's K\"unneth spectral sequence in the left half-plane.
Under appropriate circumstances (elaborated on in \cite{Sna72} and \cite{McL79}), this spectral sequence converges to $K_{G}^{*}(X\times Y)$.
In particular, the $E_{2}$ page contains the tensor product
$$\Tor^{0}_{R(G)}\bigl(K_{G}^{*}(X),K_{G}^{*}(Y)\bigr)
= K_{G}^{*}(X) \otimes_{R(G)} K_{G}^{*}(Y)$$
over the representation ring. 

In \cite[\S 9]{Hod75}, Hodgkin considers the special case $Y = G$, where
\begin{align*}
	&K_{G}^{*}(G) \cong K^{*}(\pt) \cong \Z \text{, and}\\
	&K_{G}^{*}(X\times G) \cong K^{*}(X).
\end{align*}
The K\"unneth spectral sequence then starts with $E_{2}(X,G)$ page
containing
$$\Tor^{0}_{R(G)}\bigl(K_{G}^{*}(X),\Z \bigr) = K_{G}^{*}(X) \otimes_{R(G)}\Z,$$
and since $Y=G$ is a free $G$-space, the spectral sequence converges to $K^{*}(X)$, provided that $G$ is compact, connected, and $\pi_{1}(G)$ is torsion-free (see \cite[Theorem 8.1]{Hod75}). The map (\ref{eq:weak-formality}) is then the edge homomorphism for this spectral sequence, so weak equivariant formality is precisely the statement that the edge homomorphism is an isomorphism.
This leads us to the following alternative spectral sequence condition sufficient to guarantee weak equivariant formality:

\begin{proposition}\label{prop:tor}
Let $G$ be a compact connected Lie group with $\pi_{1}(G)$ torsion-free.
A compact $G$-space $M$ is weakly equivariant formal for $K$-theory if
$$\Tor^{p}_{R(G)}\bigl( K_{G}^{*}(M), \Z\bigr) = 0$$
for all $p\neq 0$, i.e., the higher $R(G)$-torsion vanishes.
\end{proposition}

We will first use this criterion to establish weak equivariant formality for Hamiltonian $T$-spaces, where $T$ is a torus. Then, to make the transition from abelian to non-abelian actions,
we will use the following lemma (see also \cite{MM86}):

\begin{lemma}\label{lemma:g-to-t}
	Let $G$ be a compact connected Lie group with $\pi_{1}(G)$ torsion-free,
	and let $M$ be a compact $G$-space. Then $M$ is weakly equivariantly formal
	with respect to $G$ if and only if it is weakly equivariantly formal with
	respect to a maximal torus $T\subset G$.
\end{lemma}

\begin{proof}
We begin by observing that
\begin{equation}\label{eq:g-mod-t-times-m}
K^{*}_{T}(M) \cong K^{*}_{G}(G/T\times M).
\end{equation}
At the level of
vector bundles, a $T$-equivariant bundle $\pi: E \to M$ induces a
$G$-equivariant bundle $G\times_{T}E$ over $G/T\times M$. Here $G$ acts on the left on the $G$-factor of $G \times_{T}E$, while it acts via the diagonal action on $G/T \times M$. The projection map is given by $(g,v) \mapsto \left(gT,\,g\cdot \pi(v)\right)$
for $g\in G$ and $v\in E$, which is well defined since for any $t\in T$,
$$\left(gt^{-1},\,tv\right) \sim (g,v) \longmapsto \left(gt^{-1}T,\,gt^{-1}t\cdot \pi(v)\right) = \left(gT,\,g\cdot \pi(v)\right),$$
and is $G$-equivariant since for any $h\in G$,
$$h\cdot(g,v) = (hg,v) \longmapsto \left(hgT,\,hg\cdot\pi(v)\right) = h\cdot\left(g,\,g\cdot\pi(v)\right).$$
(Note that this bundle has the same total space as in Lemma~\ref{lemma:bundle}, but different base and projection.)
Conversely, given a $G$-equivariant bundle over $G/T \times M$, its restriction
to the identity coset $eT \times M$ is a $T$-equivariant bundle over $M$. These maps of vector bundles are inverses of each other and extend to $K$-theory.

In \cite[\S 3]{Sna72}, Snaith argues that under a technical hypothesis, later verified by McLeod in \cite{McL79}, Hodgkin's K\"unneth spectral sequence yields the identity
$$K^{*}_{G}(G/T \times M) \cong R(T) \otimes_{R(G)} K^{*}_{G}(M),$$
i.e., in this case the higher $\Tor$ groups vanish and the K\"unneth spectral sequence collapses.
Combining this with the identity (\ref{eq:g-mod-t-times-m}), we find that $G$-equivariant and $T$-equivariant $K$-theories of $M$ are related by
$$K^{*}_{T}(M) \cong K^{*}_{G}(M)\otimes_{R(G)} R(T).$$
Comparing this to our weak equivariant formality condition \eqref{eq:weak-formality}, we obtain
$$ K_{T}^{*}(M) \otimes_{R(T)} \Z
	\cong \left( K_{G}^{*}(M) \otimes_{R(G)} R(T) \right) \otimes_{R(T)} \Z
	\cong K_{G}^{*}(M) \otimes_{R(G)} \Z,$$
and thus weak equivariant formality with respect to $G$ and $T$ are equivalent.
\end{proof}

As an immediate corollary to this lemma, when considering equivariant formality for $K$-theory, we can replace the Lie group $G$ with any compact connected subgroup of maximal rank in $G$. We are now ready to prove our equivariant formality theorem. The basic idea of the proof, using the fact that a generic component of the moment map is equivariantly perfect with respect to a maximal torus, is based on \cite[Proposition 5.8]{Kir84}. However, we replace Kirwan's argument involving Poincar\'{e} polynomials with an argument examining the higher $\Tor$ groups. There is a similar proof in the context of algebraic $K$-theory in \cite[\S 5]{VV03}, using the Bialynicki-Birula stratification and replacing Hodgkin's topological K\"unneth spectral sequence with an algebraic $K$-theory version due to Merkur'ev \cite{Mer98}.

\begin{theorem}\label{theorem:formality}
Let $G$ be a compact connected Lie group with $\pi_{1}(G)$ torsion-free. If $M$ is a compact Hamiltonian $G$-space then $M$ is weakly equivariantly formal for $K$-theory, the forgetful map $K_{G}^{*}(M)\to K^{*}(M)$ is surjective, and every complex vector bundle over $M$ admits a stable equivariant structure.
\end{theorem}

\begin{proof}
Let $M$ be a compact Hamiltonian $G$-space. By Lemma~\ref{lemma:g-to-t} we need only establish weak equivariant formality with respect to a maximal torus $T\subset G$. Choosing a generic element $X\in\mathfrak{t}\subset\g$ in the corresponding Cartan subalgebra, the component $\mu^{X} = \langle \mu,X \rangle$ of the moment map is a Morse-Bott function whose critical sets are the connected components of the fixed point set $M^{T}$ (see \cite[Lemma 2.2]{Ati82}). Since the normal bundle to the fixed point set automatically satisfies the conditions of the Atiyah-Bott lemma (see Remark~\ref{remark:normal}), we find that the function $\mu^{X}$ is equivariantly perfect for $K$-theory. We can therefore compute $K^{*}_{T}(M)$ as was done using formula (\ref{eq:short-exact-sequence-c}) by taking a series of extensions of the form
$$\xymatrix{
0 \ar[r] & K_{T}^{*}\bigl(M^{T}_{i}\bigr) \ar[r]  & 
K_{T}^{*}\bigl(M_{\leq i}\bigr) \ar[r] & 
K_{T}^{*}\bigl(M_{<i}\bigr) \ar[r] & 0},
$$
where $i$ indexes the connected components of the fixed point set $M^{T}$. Since $T$ acts trivially on its fixed point set, for each of the critical sets we have
$$K_{T}^{*}\bigl(M^{T}_{i}\bigr) \cong K^{*}\bigl(M^{T}_{i}\bigr) \otimes R(T),$$
and thus their higher $\Tor$ groups vanish, i.e.,
$$\Tor^{p}_{R(T)}\bigl(K_{T}^{*}(M^{T}_{i}),\Z\bigr) = 0 \text{ for } p \neq 0.$$
By the long exact sequence for the $\Tor$ functor (see \cite[\S VIII.1]{CE56} for a discussion of the algebra), this property continues to hold upon taking a series of extensions, and so for all of $M$ we have 
$$\Tor_{R(T)}^{p}\bigl(K_{T}^{*}(M),\Z\bigr) = 0 \text{ for } p\neq 0$$
It follows from Proposition~\ref{prop:tor} that $M$ is weakly equivariantly formal for $K$-theory. The remaining assertions follow from Proposition~\ref{prop:formality}.
\end{proof}

%Applying Proposition~\ref{prop:formality}, we obtain

%\begin{corollary}\label{cor:formality}
%Let $M$ be a compact Hamiltonian $G$-space with respect to a compact connected Lie group $G$ with $\pi_{1}(G)$ torsion-free. Then the forgetful map $K_{G}^{*}(M)\to K^{*}(M)$ is surjective and every complex vector bundle over $M$ admits a stable equivariant structure.
%\end{corollary}

\begin{remark}\label{rem:kirwan-ginzburg}
Kirwan proves the standard version of equivariant formality for rational cohomology in
\cite[\S 5]{Kir84}. The Leray-Serre spectral sequence gives an upper
bound on the size (in terms of the Betti numbers) of the equivariant
cohomology $H_{G}^{*}(M;\Q)$. On the other hand, knowing that the
components of the moment map are equivariantly perfect Morse-Bott
functions, the Morse inequalities for the non-equivariant cohomology
provide a lower bound on the size of $H_{G}^{*}(M;\Q)$, establishing
the result. As a corollary, Kirwan notes that the components of the
moment map are also perfect in the non-equivariant case. Ginzburg in
\cite{Gin87} works with de Rham cohomology and reverses this argument,
showing that the moment map components are perfect and then using that
to establish equivariant formality. As a consequence, on a compact
Hamiltonian $G$-space, any rational cohomology class admits an
equivariant extension, and any closed form admits an equivariantly
closed extension.
\end{remark}

A deeper question is whether we can remove the word \emph{stable} from Theorem~\ref{theorem:formality}. In other words, does every complex vector bundle over a Hamiltonian $G$-space admit a lift of the $G$-action, without requiring any additional degrees of freedom? In \cite{Kos70}, Kostant shows that the moment map contains precisely the data required to lift the $G$-action on $M$ to a $G$-action on the prequantum line bundle $L\to M$ satisfying $c_{1}(L) = [\omega]$.
%On a Hamiltonian $G$-space, we can lift the $G$-action to the prequantization line bundle $L\to M$ with $c_{1}(L) = [\omega]$; in
%fact, the moment map contains precisely the data required to lift the
%$G$-action on $M$ to a $G$-action on $L$ (see \cite{Kos70}).
%However,
%we would like to show that such a lift exists for all vector bundles.
%Since the Chern character is an isomorphism from rational $K$-theory to
%rational cohomology, where we have equivariant formality, the
%map $K_{G}(M) \otimes \Q \to K(M) \otimes \Q$ must be
%surjective. Decomposing $K(M) = F \oplus T$ into free and torsion
%components respectively, tensoring with the rationals eliminates the
%torsion, giving us $K(M) \otimes \Q = F \otimes \Q$. For any
%non-equivariant bundle $E$ over $M$ corresponding to a free class $[E]
%\in F$, it follows that some multiple $k\,[E]$ for $k\geq 1$ lies in
%the image of $K_{G}(M)$, and thus $E^{\oplus k}$ admits an equivariant
%structure. However, passing to rational cohomology tells us nothing
%about torsion bundles. On the other hand, for the case of line bundles
%we can tackle the torsion directly, and
In the following we show that
%prove the partial result that 
every complex line bundle admits an
equivariant structure, at least with respect to a finite cover of $G$. (See
\cite{HY76, MIR01} for related discussions of lifts of group actions
to line bundles.)
% GDL 3/11/05

\begin{theorem}\label{equivariant-line-bundles}
Let $G$ be a compact connected Lie group with $\pi_{1}(G)$ torsion-free, and let $M$ be a compact Hamiltonian $G$-space.
\begin{enumerate}
\item
Every complex line bundle over $M$ admits a fiberwise linear lift of the $G$-action on $M$.
\item
If $M$ is connected, then the isomorphism classes of the equivariant structures on a fixed complex line bundle over $M$ are in one-to-one correspondence with the character group of $G$.
\item
If $G$ is semi-simple, then each complex line bundle over $M$ has a unique equivariant structure.
\end{enumerate}
\end{theorem}

\begin{proof}
We recall that the Picard group $\Pic(M)$ of isomorphism classes of complex line bundles over $M$ is a retract of $K(M)$. In particular, the determinant homomorphism $\det : K(M) \twoheadrightarrow \Pic_(M)$
induced by the map taking a complex vector bundle $E$ to its top exterior power $\det E = \Lambda^{\rank E}(E)$ is surjective. The same holds in the equivariant case, and we have a retraction $\det_{G} : K_{G}(M) \twoheadrightarrow \Pic_{G}(M)$. 
Using Theorem~\ref{theorem:formality}, we can construct the commutative diagram
$$\xymatrix{
	K_{G}(M) \ar@{->>}[r] \ar@{->>}[d]_{\det_{G}} & K(M) \ar@{->>}[d]^{\det} \\
	\Pic_{G}(M) \ar[r] & \Pic(M)
}$$
from which we see that the forgetful map $\Pic_{G}(M) \to \Pic(M)$ is also surjective.

The kernel of the forgetful map consists of all equivariant structures on the trivial complex line bundle $M \times \C$. Restricted to a single $G$-orbit $G \cdot x$ in $M$, such $G$-actions are necessarily of the form $(G \cdot x) \times \C_{\lambda}$, where $\C_{\lambda}$ is a 1-dimensional representation of $G$.  Over all of $M$, an equivariant structure on the trivial bundle corresponds to a continuous map
$M/G \to \Pic_{G}(\pt)$ from the orbit space to the character group of all such 1-dimensional representations (i.e., continuous homomorphisms from $G$ to the circle group). Since the character group $\Pic_{G}(\pt)$ is discrete, the representation $\C_{\lambda}$ must be constant on each connected component of $M$. The equivariant structures therefore correspond to $\Pic_{G}(\pi_{0}M)$, and if $M$ is connected they correspond to the character group $\Pic_{G}(\pt)$.

In summary, we obtain the short exact sequence of Picard groups
$$1 \longrightarrow \Pic_{G}(\pi_{0}M) \longrightarrow \Pic_{G}(M) \longrightarrow \Pic(M) \longrightarrow 1$$
whenever $M$ is a compact Hamiltonian $G$-space with respect to a compact connected Lie group $G$ with $\pi_{1}(G)$ torsion-free. For the final assertion, we recall that a compact semi-simple Lie group has no one-dimensional representations other than the trivial representation, and thus its character group is trivial. 
\end{proof}

\begin{example}
%Note that we required that $M$ be symplectic and that the $G$-action be Hamiltonian only to prove surjectivity onto the free elements. This argument shows that we always have surjectivity onto the torsion elements, regardless of the symplectic structure or existence of a moment map, provided that $G$ is connected and has no torsion in its fundamental group. On the other hand,
If the $G$-action on $M$ is not Hamiltonian, then equivariant formality may fail.
% for the free elements.
Let $M = T^{2}$ with
$G = S^{1}$ acting by rotating one of the two circles freely. This is not a Hamiltonian action. The equivariant $K$-theory is
$$K^{*}_{S^{1}}(T^{2}) \cong K^{*}(T^{2}/S^{1}) \cong K^{*}(S^{1}),$$
which does not surject onto the standard $K$-theory $K^{*}(T^{2})$.
Also, since $S^{1}$ acts freely on $T^{2}$, the representation ring $R(S^{1})$acts on $K_{S^{1}}(T^{2})$ via its image under the augmentation homomorphism, and we obtain
$$K^{*}_{S^{1}}(T^{2}) \otimes_{R(S^{1})} \Z \cong K_{S^{1}}^{*}(T^{2}).$$
For line bundles we have
$$\Pic_{S^{1}}(T^{2}) \cong \Pic(T^{2}/S^{1}) \cong \Pic(S^{1}) = 0,$$
while $\Pic(T^{2})\cong \Z$. Thus there exist non-trivial line bundles
over $T^{2}$ which do not admit equivariant structures, as all equivariant line bundles must be trivial.
\end{example}

%\section{Examples}\label{sec:examples}

%\subsection{Coadjoint orbits}\label{subsec:coadjOrbits}

%Maybe a computation with a $U(n,\C)$-coadjoint orbit weight variety? ?? 

%\subsection{Hypertoric varieties}\label{subsec:hypertorics}

%It might be interesting to look at hypertorics as a \(T\times
%S^1\)-Hamiltonian space using the ``first'' symplectic structure, and
%looking at the {\em symplectic} (NOT hyperKahler) reduction. Here you
%don't just get a point -- but I'm not sure what it is that you get. ..

\begin{remark}
The reader may be curious why all the results in this section require the fundamental group of $G$ to be free abelian. This condition appears in Hodgkin's
derivation of the K\"unneth spectral sequence in \cite{Hod75}, and in Remark~\ref{rem:spectral-sequence} below we show explicitly how this condition is used
in one particular case. Most significantly, the condition on the torsion of $\pi_{1}(G)$ is the $K$-theoretic analogue of the condition that $H^{*}(G;\Z)$ be torsion free which we encountered at the end of Section~\ref{sec:AtiyahBott} when discussing the torsion constraints of the cohomological Atiyah-Bott lemma. In fact, by Hodgkin's Theorem (see \cite{Hod67,McL79}), if a compact connected Lie group has no torsion in its fundamental group, then its $K$-theory is automatically torsion-free.

We observe that if $\pi_{1}(G)$ contains torsion,
then by the Hurewicz theorem the homology group $H_{1}(G;\Z)$ also contains the same torsion.
By the universal coefficient theorem, we then obtain torsion in $H^{2}(G;\Z)$,
which in turn transgresses to give torsion in $H^{3}(BG;\Z) \cong H_{G}^{3}(\pt)$.
However, twistings of $K$-theory, as described in \cite{AS03,FHT02} are classified by elements of $H_{G}^{3}(M)$. Our condition on the torsion of $\pi_{1}(G)$ ensures that $H_{G}^{3}(\pt)$ vanishes, and thus we do not encounter any twistings of $K_{G}^{*}(M)$ arising solely from the Lie group $G$. We expect that many of our results of this section still hold without the torsion constraint, albeit in terms of twisted $K$-theory.

\end{remark}

\appendix

\section{Equivariant Cohomology and Line Bundles}\label{appendix-line-bundle}

Our results from Section~\ref{formality} involving equivariant line bundles and the equivariant Picard group can also be stated in terms of equivariant cohomology. We recall that $G$-equivariant line bundles over a compact $G$-space $M$ are classified up to isomorphism by $H^{2}_{G}(M;\Z)$. This is a non-trivial result, which appears in the literature in several sources, including \cite[Theorem C.47]{GGK} and \cite{LMS83}. Here we present an elementary
% ``dead elementary''\footnote{according to J. Frank Adams when it was presented to him on a napkin.}
proof due to Peter Landweber. We recall that for compact $M$, the first Chern class $c_{1} : \Pic(M) \to H^{2}(M;\Z)$ is an isomorphism taking tensor products of line bundles to sums of cohomology classes, as we have
$$\Pic(M) \cong [M,\BU_{1}] \cong [M,K(\Z,2)] \cong H^{2}(M;\Z).$$
We now establish the equivariant version of this statement.

\begin{theorem}\label{equivariant-chern-class}
	Let $G$ be a compact Lie group acting on a compact manifold $M$.
	The equivariant first Chern class
	\begin{equation}\label{eq:equivariant-c1}
	c_{1}^{G} : \Pic_{G}(M) \to H^{2}_{G}(M;\Z)
	\end{equation}
	is an isomorphism between the equivariant Picard group of isomorphism
	classes of complex $G$-line bundles on $M$ and the equivariant cohomology
	$H^{2}_{G}(M;\Z) = H^{2}(M_{G};\Z)$,
	defined in terms of the Borel construction $M_{G} = EG \times_{G} M$.
\end{theorem}

\begin{proof}
To prove this result, we compare the $G$-equivariant first Chern class (\ref{eq:equivariant-c1}) on $M$ to the standard first Chern class for the Borel construction $M_{G}$.
We recall that the $K$-theory of a general (possibly non-compact) space $X$ is defined homotopically by
$$K(X) := [X,\,\Z \times BU].$$
For any $X$, the first Chern class gives us a retraction
\begin{equation}\label{eq:general-c1}
c_{1} : K(X) \to H^{2}(X;\Z),
\end{equation}
taking the homotopy class $[f]$ of a map $f : X \to 
\Z \times \BU$ to the pullback $f^{*}a$ of the preferred generator $a\in H^{2}(\BU;\Z)$. (Alternatively, we
can take the composition of $f$ with a representative of the class $a\in[\BU,K(\Z,2)]$, whose homotopy class gives $c_{1}[f]\in[X,K(\Z,2)]$.)
The right inverse of the homomorphism (\ref{eq:general-c1}) is the injection
\begin{equation}\label{eq:retraction-inverse}
H^{2}(X;\Z) \cong [X,\BU_{1}] \longrightarrow [X,\,\Z \times \BU] \cong K(X),
\end{equation}
induced by the inclusion $U_{1} \hookrightarrow U$, with $\BU_{1}$ mapping into the $1 \times \BU$ connected component.
We observe that if $X$ is compact, then this homomorphism (\ref{eq:retraction-inverse}) is the composition
$$H^{2}(X;\Z) \xrightarrow{c_{1}^{-1}} \Pic(X) \longrightarrow K(X)$$
of the inverse of the first Chern class with the homomorphism taking a complex line bundle $L$ to its $K$-theory class $[L]$. However, we shall consider these maps in the non-compact case $X= M_{G}$. On the other hand, the equivariant Picard group $\Pic_{G}(M)$ is a retract of the equivariant $K$-theory $K_{G}(M)$, as we discussed in the proof of Theorem~\ref{equivariant-line-bundles}.%
%mentioned, we have an inclusion $$\Pic_{G}(M) \hookrightarrow
%K_{G}(M)$$ taking $L \mapsto [L]$. Its left inverse is the determinant
%homomorphism $$\det : K_{G}(M) \twoheadrightarrow \Pic_{G}(M),$$
%induced by the map taking a complex $G$-bundle $E$ to its top exterior
%power $\det E = \Lambda^{\mathrm{rank}\,E}(E)$.
\footnote{The map \(c_1 :
K(M_G) \to H_G^2(X;\Z)\) 
can also be interpreted geometrically as the map \(\det :
K(M_G) \twoheadrightarrow \Pic(M_G).\)}

We recall that the equivariant Chern classes of a $G$-bundles $E \to M$ are defined by %in terms of the Borel construction by
$$c_{i}^{G}(E) := c_{i}( EG \times_{G} E ) \in H^{2i}(M_{G}) = H^{2i}_{G}(M).$$
In terms of $K$-theory, let $\alpha : K_{G}(M) \to K(M_{G})$ denote the Atiyah-Segal homomorphism,
induced by the map taking a $G$-bundle $E$ to the associated
bundle $EG\times_{G}E $ over $M_{G} = EG\times_{G} M$. We then have
$c^{G}_{i}[E] = c_{i}(\alpha[E])$, which gives us the following commutative diagram:
%\begin{equation*}
%\xymatrix{
%	\Pic_{G}(M) \ar[r]^{c_{1}^{G}} \ar@{^{(}->}[d] & H^{2}_{G}(M;\Z) \\
%	K_{G}(M)  \ar[r]^{\alpha} & K(M_{G}) \ar@{->>}[u]_{c_{1}}
%}
%\end{equation*}
\begin{equation*}
\xymatrix{
	\Pic_{G}(M) \ar[r]^{c_{1}^{G}} \ar@{^{(}->}[d]
		& H^{2}_{G}(M;\Z) \ar@{^{(}->}[d] \\
	K_{G}(M)  \ar[r]^{\alpha} \ar@{->>}[d]_{\det}
		& K(M_{G}) \ar@{->>}[d]^{c_{1}} \\
	\Pic_{G}(M) \ar[r]^{c_{1}^{G}}
		& H^{2}_{G}(M;\Z)
}
\end{equation*}
where both vertical compositions are identity maps. From this we observe that the equivariant first Chern class $c_{1}^{G}$ is a retract of the Atiyah-Segal homomorphism $\alpha$.

In Lemma~\ref{determinant} below, we prove that the determinant homomorphism
extends to the completion $K_{G}(M)^{\wedge}$ of the equivariant $K$-theory $K_{G}(M)$ with respect to the augmentation ideal $I(G)\cdot K_{G}(M)$. It follows that $\Pic_{G}(M)$ is a retract of $K_{G}(M)^{\wedge}$ as well.
Our commutative diagram then becomes:
\begin{equation*}
\xymatrix{
	\Pic_{G}(M) \ar[rr]^{c_{1}^{G}} \ar@{^{(}->}[d]
		&& H^{2}_{G}(M;\Z) \ar@{^{(}-->}[d] \\
	K_{G}(M)  \ar@{^{(}->}[r] \ar@{->>}[d]_{\det}
		& K_{G}(M)^{\wedge} \ar@{-->}[r]^{\widehat{\alpha}}_{\cong} \ar@{-->>}[dl]_{\widehat{\det}} & K(M_{G}) \ar@{->>}[d]^{c_{1}} \\
	\Pic_{G}(M) \ar[rr]^{c_{1}^{G}}
		&& H^{2}_{G}(M;\Z)
}
\end{equation*}
%\begin{equation*}
%\xymatrix{
%	\Pic_{G}(M) \ar[r]^{c_{1}^{G}} \ar@{^{(}->}[d] & H^{2}_{G}(M;\Z) \ar@<-1ex>@{_{(}-->}[d] \\
%	K_{G}(M) \ar[r]^{\alpha} & K(M_{G}) \ar@{-->>}[ul]_{\widehat{\det}} \ar@{->>}@<-1ex>[u]_{c_{1}}
%}
%\end{equation*}
The equivariant first Chern class $c_{1}^{G}$ is therefore a retract of the map $\widehat{\alpha}: K_{G}(M)^{\wedge} \to K(M_{G})$, which by the Atiyah-Segal completion theorem is an isomorphism (see \cite{Ati-Seg}). As a consequence, the map $c_{1}^{G}$ is an isomorphism, with inverse given by the composition of the dashed arrows.
\end{proof}

To complete the proof, we must establish the following technical lemma, whose proof requires only linear algebra. % and the equivariant splitting principle.

\begin{lemma}\label{determinant}
The determinant homomorphism $\det : K_{G}(M) \to \Pic_{G}(M)$ extends to a homomorphism $\widehat{\det} : K_{G}(M)^{\wedge} \to \Pic_{G}(M)$ on the completion at the augmentation ideal.
\end{lemma}

\begin{proof}
Consider the augmentation homomorphism $K_{G} \to \Z$ taking a virtual $G$-module to its integral dimension. The augmentation ideal $I(G)$ is the kernel of this map. The completion of the equivariant $K$-theory with respect to the $I(G)$-adic topology is given by
$$K_{G}(M)^{\wedge} = \lim_{\longleftarrow} K_{G}(M) \,/\, I(G)^{n}\cdot K_{G}(M).$$
We claim that the determinant is trivial on $I(G)^{2}\cdot K_{G}(M)$, and it thus extends to the $I(G)$-adic completion $K_{G}(M)^{\wedge}$.
 Let $E \cong P_{1} \oplus \cdots \oplus P_{p}$ and $F \cong Q_{1} \oplus \cdots \oplus Q_{q}$ be two 
complex $G$-bundles, each splitting as a direct sum of complex $G$-line bundles.
Their tensor product is then $E \otimes F \cong \bigoplus_{i,j} P_{i} \otimes Q_{j}$, which has determinant
\begin{equation}\label{eq:det-tensor}
\det(E \otimes F) = \Bigl( {\prod}_{i} P_{i}^{\otimes q} \Bigr) \otimes \Bigl( {\prod}_{j} Q_{j}^{\otimes p} \Bigr)
= (\det E)^{\mathrm{rank}\,F} \otimes (\det F)^{\mathrm{rank}\,E}.
\end{equation}
By the equivariant splitting principle (see \cite{Seg68}), the formula (\ref{eq:det-tensor}) holds for all $G$-bundles. Furthermore, since the determinant of a sum is the product of the determinants in the Picard group, we can extend the formula (\ref{eq:det-tensor}) to virtual $G$-bundles in place of $E$ and $F$.
In particular
if $V, W$ are both virtual $G$-modules with rank $0$, and $E$ is an arbitrary virtual $G$-bundle, then $W \otimes E$ has rank $0$, and we obtain
$$\det(V \otimes W \otimes E) = (\det V)^{0} \otimes (\det W \otimes E)^{0} =1,$$
the identity element in $\Pic_{G}(M)$.
\end{proof}

%We note that the proof of Lemma~\ref{determinant} does not actually use
%the equivariant structure, other than to invoke the equivariant splitting principle. Without the equivariance, Lemma~\ref{determinant} becomes a statement about the completion of $K(M)$ with respect to the trivial ideal,
%% GDL 3/11/05
%% $\tilde{K}(M)$ (assuming that $M$ is connected). However, since the reduced $K$-theory is nilpotent, taking such a completion leaves the $K$-theory unchanged,
%and the result is vacuous.

\begin{theorem}\label{theorem:h2}
Let $G$ be a compact connected Lie group with $\pi_{1}(G)$ torsion-free, and let $M$ be a compact Hamiltonian $G$-space. Then we have a short exact sequence
\begin{equation}\label{eq:h2-sequence}
0 \longrightarrow H_{G}^{2}(\pt)\otimes H^{0}(M;\Z) \longrightarrow
H_{G}^{2}( M;\Z) \longrightarrow H^{2}(M;\Z) \longrightarrow 0,
\end{equation}
i.e., the forgetful map is surjective and its kernel is isomorphic to a direct sum of one copy of $H_{G}^{2}(BG;\Z)$ for each connected component of $M$. From this, we obtain the following:
\begin{enumerate}
\item
Every cohomology class in $H^{2}(M;\Z)$ admits an equivariant extension.
\item
If $M$ is connected, then the isomorphism classes of the equivariant structures on a given cohomology class in $H^{2}(M;\Z)$ are in one-to-one correspondence with $H^{2}_{G}(\pt;\Z)$.
\item
If $G$ is semi-simple, then the forgetful map is an isomorphism, and thus every class in $H^{2}(M;\Z)$ admits a unique equivariant extension.
\end{enumerate}
\end{theorem}
\begin{proof}
These results follows immediately from Theorem~\ref{equivariant-line-bundles}, together with the first Chern class isomorphisms $\Pic(M) \cong H^{2}(M;\Z)$ and $\Pic_{G}(M) \cong H_{G}^{2}(M;\Z)$ from Theorem~\ref{equivariant-chern-class}.
\end{proof}

\begin{remark}\label{rem:spectral-sequence}
Alternatively, we can prove Theorem~\ref{theorem:h2} using the Leray-Serre spectral sequence for the fibration $M \to M_{G} \to BG$. As we noted in Remark~\ref{rem:kirwan-ginzburg}, a Hamiltonian $G$-space $M$ is equivariantly formal for
rational cohomology, so with rational coefficients this spectral sequence collapses at the $E_{2}$ term, and the rational differentials
$d_{k}^{\Q}$ all vanish. In particular, the image of the forgetful
map $H_{G}^{2}(M) \to H^{2}(M)$ is $E_{\infty}^{0,2}$, and
thus this map is surjective if and only if $H^{2}(M)$ persists to
$E_{\infty}$, or equivalently if all differentials $d_{k}$ vanish on
$H^{2}(M)$. Using integral cohomology, since $\pi_{1}(G)$ is
free, we see that $H_{1}(G;\Z)$ is free, and that
$H^{2}(G;\Z)$ is likewise free by the universal coefficient
theorem. Since $G$ is connected, its classifying space $BG$ is simply
connected, and it follows that $H^{1}(BG;\Z) = H^{3}(BG;\Z) = 0$ and
$H^{2}(BG;\Z)$ is free. Applying the universal coefficient theorem
once again, we see that the cohomology $H^{1}(M;\Z)$ of the fiber
is also free.  The $E_{2}$ page is then
%For a class in $H^{2}(M;\Z)$ to persist to the equivariant cohomology $H^{2}_{G}(M;\Z)$, it must survive
%the $E_{2}$ and $E_{3}$ pages of the spectral sequence.  At the $E_{2}$ page, we have % the differential
%%$$d_{2} : E_{2}^{0,2} = H^{2}(M;\Z) \longrightarrow E_{2}^{2,1} \cong H^{1}(M;\Z) \otimes H^{2}(BG;\Z).$$
% GDL 3/11/05
\begin{equation*}
\xymatrix{
	H^{2}(M;\Z) \ar@<1ex>[drr]^{d_{2}} \ar[ddrrr]_(.7){d_{3}} \\
	H^{1}(M;\Z) & 0 & H^{2}\bigl(BG;\,H^{1}(M;\Z)\bigr) & 0 \\
	H^{0}(M;\Z) & 0 & H^{2}\bigl(BG;\,H^{0}(M;\Z)\bigr) & 0
}
\end{equation*}
where the transgression $d_{3}$ actually acts on the kernel of $d_{2}$.
Since the image of $d_{2}$ is free, and since tensoring with the rationals gives
$(\mathrm{Im}\,d_{2}) \otimes \Q = \mathrm{Im}\,d_{2}^{\Q} = 0,$
we see that $d_{2}$ vanishes on $H^{2}(M;\Z)$. Furthermore, the transgression $d_{3}$ and $d_{k}$ for $k > 3$ all clearly vanish on $H^{2}(M;\Z)$.
%However, the image of $d_{2}$ is free, and thus it vanishes on the torsion, while the transgression $d_{3}$ clearly vanishes everywhere.
% At the $E_{3}$ page, since $H^{3}(BG) = 0$, the transgression $d_{3} : E_{3}^{0,2} \to E_{3}^{3,0}$ must vanish.
As a result, all of $H^{2}(M;\Z)$ persists to $E_{\infty}$, and thus $H_{G}^{2}(M;\Z)$ surjects onto $H^{2}(M;\Z)$.

By a similar argument, the transgression $d_{2} : H^{1}(M;\Z) \to H^{2}(BG;H^{0}(M;\Z))$ vanishes as its range is free, and so $H^{2}(BG;H^{0}(M;\Z))$ persists to $E_{\infty}$. We thus obtain the short exact sequence
$$ 0 \longrightarrow H^{2}\bigl(BG;H^{0}(M;\Z)\bigr) \longrightarrow H^{2}(M_{G};\Z) \longrightarrow H^{2}(M;\Z) \longrightarrow 0.$$
%As a consequence, we see that $H^{2}_{G}(\mathrm{pt};H^{0}(M;\Z))$ is the kernel
%of the restriction map $H_{G}^{2}(M;\Z) \to H^{2}(M;\Z)$, and thus the
%possible equivariant extensions of a class in $H^{2}(M;\Z)$ correspond
%to $H^{2}_{G}(\mathrm{pt};H^{0}(M;\Z))$.
which is equivalent to (\ref{eq:h2-sequence}). We note that this argument shows explicitly why we require that $G$ have a free fundamental group, as well as giving a different formulation of the kernel of the forgetful map.
\end{remark}

\section{Morse-Kirwan functions}\label{morse-kirwan}

In \cite[\S 10]{Kir84}, Kirwan extends the results of Morse theory (or more properly, Morse-Bott theory) to a larger class of Morse functions, which she calls ``minimally degenerate'' and which are now commonly called Morse-Kirwan functions (or sometimes ``Morse in the sense of Kirwan'').

\begin{definition}
A \emph{Morse-Kirwan function} on a manifold $M$ is a smooth function $f : M \to \R$ whose critical set decomposes as the union of disjoint closed critical subsets on which $f$ takes constant values, such that each of these critical subsets $C$ is contained in a locally closed submanifold $\Sigma_{C}$
with the following properties:
\begin{enumerate}
\item The critical subset $C$ is the subset of $\Sigma_{C}$ minimizing the function $f$.
\item The neighborhood $\Sigma_{C}$ has an orientable normal bundle.
\item The tangent space $T_{x}\Sigma_{C}$ to $\Sigma_{C}$ at $x\in C$ is maximal among all subspaces of the tangent space $T_{x}M$ on which the Hessian $H_{x}(f)$ is positive semidefinite.
\end{enumerate}
\end{definition}

For such a Morse-Kirwan function, the critical sets may not be manifolds, as they are for Morse-Bott functions, but rather they may have singularities. However, we are allowed to resolve these singularities by flowing upwards slightly to obtain a manifold neighborhood of the critical set, whose normal bundle consists of the downward directions needed for Morse theory.

Kirwan shows in \cite[\S 10]{Kir84} that such minimally degenerate Morse functions yield smooth stratifications, and thus a form of Morse theory still holds, at least at the level of the stratification. The only complication is that for a general Morse-Kirwan function,
%uses such a Morse function, together with a Riemannian metric, to obtain a negative gradient flow on the manifold $M$. She then
%constructs a smooth stratification $M = \bigsqcup S_{C}$ indexed by the connected components $\{C\}$ of the critical set of $f$, where each stratum $S_{C}$ is the set of all points $x\in X$ which flow downward to the critical subset $C$ along their paths of steepest descent. (
%Note that for a Morse-Kirwan function,
the negative gradient flow of a point need not converge to unique
limit; rather, its set of limit points is a closed, connected subset
of the critical set of $f$. (On the other hand, the gradient flow of
the norm square of the moment map does indeed converge, as described
in unpublished work of Duistermaat, reproduced in \cite{Ler04}, and
also proved in \cite[Appendix B]{Woo02}.)
Instead, Kirwan argues that the inclusion $C \hookrightarrow S_{C}$ of a critical set into its stratum is a retraction in (equivariant) \v{C}ech cohomology. We now show, reprising carefully Kirwan's argument from \cite[\S 10]{Kir84}, that the inclusion $C \hookrightarrow S_{C}$ is also a retraction in (equivariant) $K$-theory. Let $M$ be a $G$-manifold, and consider a $G$-invariant Morse-Kirwan function $f : M \to \R$.  For small $\delta \geq 0$, we consider the compact neighborhood
$$N_{\delta} = \{ x\in S_{C} \,|\, f(x) \leq f(C) + \delta \}$$
of $C$ in $S_{C}$. Note that $N_{0} = C$, and in the following we approximate the critical set $C$ as $\lim_{\delta\to 0}N_{\delta}$.

\begin{lemma}\label{lemma:retraction}
	The gradient flow of $f$ induces a $G$-equivariant deformation retraction
	of the stratum $S_{C}$ onto each of the neighborhoods $N_{\delta}$ for $\delta > 0$.
\end{lemma}
\begin{proof}
Let $\gamma : M \times \R^{\geq 0} \to M$ denote the gradient flow of $f$, satisfying $\partial_{t}\,\gamma_{t}(x) = -(\grad f)|_{\gamma_{t}(x)}.$
Fixing $\delta > 0$, for each point $x\in S_{C}$, we define
$$\tau(x) := \inf \{ t \geq 0 \,|\, \gamma(x,t) \in N_{\delta} \}$$
to be the time at which the gradient flow of $x$ first enters the neighborhood $N_{\delta}$.
In particular, we have $\tau(x) = 0$ if and only if $x\in N_{\delta}$,
and we see that $\tau(x)$ is well defined as $S_{C}$ consists of all points which flow down to the critical set $C$. Leting $\epsilon : [0,1) \to [0,\infty)$ be the homeomorphism $\epsilon(s) = s/(1-s)$, we define a homotopy $F : S_{C} \times [0,1] \to S_{C}$ by
$$F_{s}(x) = \begin{cases}
	\gamma_{\tau(x)}(x) & \text{if $s=1$ or $\epsilon(s) \geq \tau(x)$,} \\
	\gamma_{\epsilon(s)} (x) & \text{otherwise.} \\
\end{cases}$$
Then $F_{0}(x)$ is the identity map on $S_{C}$, as $\epsilon(0) = 0$ and $\gamma_{0}(x) = x$. On the other hand, the map $F_{1}(x)$ is the retraction of the stratum $S_{C}$ onto the neighborhood $N_{\delta}$ obtained by taking each $x\in S_{C}$ to the point at which its gradient flow first enters $N_{\delta}$.

As for the $G$-equivariance, since the Morse-Kirwan function $f$ is $G$-invariant, the $G$-action maps the neighborhood $N_{\delta}$ to itself, and the gradient flow $\gamma_{t}(x)$ is $G$-equivariant.
For $x\in S_{C}$ and $g\in G$, it follows that $\tau(gx) = \tau(x)$,
and thus
$\gamma_{\tau(gx)}(gx) = g \gamma_{\tau(x)}(x)$. The homotopy $F_{s}(x)$ is therefore a $G$-equivariant deformation retraction
of $S_{C}$ onto $N_{\delta}$.
\end{proof}

\begin{lemma}\label{lemma:retraction-isomorphism}
	The inclusion $C \hookrightarrow S_{C}$ of a critical set into its stratum
	induces an isomorphism $K_{G}^{*}(S_{C}) \cong K_{G}^{*}(C)$ in equivariant
	$K$-theory.
\end{lemma}

\begin{proof}
By Lemma~\ref{lemma:retraction}, the inclusions $N_{\delta}\hookrightarrow S_{C}$ induce isomorphisms
$K_{G}(S_{C}) \cong K_{G}(N_{\delta})$ in equivariant $K$-theory for all $\delta > 0$. To compute the equivariant $K$-theory of the critical set, $K_{G}(C) = K_{G}(N_{0})$, observe
that the compact neighborhoods $N_{\delta}$ form a directed inverse system
with inverse limit $C$ as $\delta \to 0$:
$$C = \lim_{\longleftarrow}N_{\delta}.$$
In our case, we have inclusions $N_{\delta'}\hookrightarrow N_{\delta}$ for $\delta' < \delta$ and the simplified condition that $C = \bigcap_{\delta>0}N_{\delta}$.
Furthermore, it follows from Lemma~\ref{lemma:retraction} that the corresponding restriction maps $K_{G}(N_{\delta}) \to K_{G}(N_{\delta'})$ are isomorphisms. So,
we invoke the continuity axiom for equivariant $K$-theory, which asserts that for compact Hausdorff spaces, the equivariant $K$-theory of an inverse limit is the direct limit of the equivariant $K$-theories (see \cite[Proposition 2.11]{Seg68}), to obtain the isomorphism
$$K_{G}^{*}(C) = K_{G}^{*}\bigl(\lim_{\longleftarrow}N_{\delta}\bigr) = \lim_{\longrightarrow}K_{G}^{*}(N_{\delta})
\cong \lim_{\longrightarrow } K_{G}^{*}(S_{C})
 = K_{G}^{*}(S_{C})$$
induced by the inclusion $C\hookrightarrow S_{C}$.
\end{proof}

\begin{remark}
To establish the continuity axiom in equivariant $K$-theory, we recall
from \cite{Cap54} that the continuity axiom is equivalent to the
extension and reduction theorems. Let $A$ be a closed subset of a
compact Hausdorff space $X$. The extension theorem says that any
$K$-theory class on $A$ can be extended to a class on a compact
neighborhood of $A$. The reduction theorem says that if a $K$-theory
class on $X$ vanishes when restricted to $A$, then it likewise
vanishes on some compact neighborhood of $A$. However, for $K$-theory,
these results are both true at the level of vector bundles. See for
example the texts \cite[\S1.4]{Ati67} or \cite[proof of Proposition 2.9]{Hat03} for
discussions of reduction for standard $K$-theory. These results then
remain true when passing to the Grothendieck group. Segal argues along
these lines in \cite{Seg68}, in which he proves the continuity
property of equivariant $K$-theory in the case we require above where
the inverse limit is actually the intersection.
\end{remark}

\smallskip

%We therefore obtain the following Thom-Gysin sequence in equivariant $K$-theory:
%\begin{equation*}
%\xymatrix{
%\cdots \ar[r] & K_{G}^{*-d(C)}\bigl(C) \ar[r]  & 
%K_{G}^*\bigl(\bigsqcup_{C' \leq C} S_{C'}\bigr) \ar[r] & 
%K_{G}^*\bigl(\bigsqcup_{C' < C} S_{C'}\bigr) \ar[r] & \cdots }
%\end{equation*}
%again provided that the normal bundle to $\Sigma_{C}$ in $M$ admits a $\Spinc$-structure.
We expect a similar form of Morse theory to hold for general Morse-Kirwan functions in any cohomology theory which has Thom classes and satisfies the continuity axiom.
%, i.e., for which $\lim^{1}$ vanishes. 
In addition, when working with the norm square of the moment map, we can drop the continuity condition.

%\todo{Reference for $\lim^{1}$.}

%More precisely, since isomorphism classes of rank $k$ complex bundles are classified by the \v{C}ech cohomology group $H^{1}(X;\GL_{k}\C)$, and since \v{C}ech cohomology itself satisfies the continuity axiom, we see that the extension and reduction theorems hold for vector bundles. For our purposes, we are actually interested in $G$-equivariant $K$-theory. However, the continuity axiom likewise holds for the equivariant \v{C}ech cohomology $H^{1}_{G}(X;\GL_{k}\C) = H^{1}(X_{G};\GL_{k}\C)$, using the Borel construction $X_{G} = X \times_{G} EG$, provided that we
%approximate $EG$ as the direct limit of compact manifolds. Thus the $K$-theory of $G$-equivariant vector bundles satisfies the continuity axiom as well.

%\bibliographystyle{habbrv}
%\bibliography{ref}

\def\cprime{$'$}

\end{spacing}
\end{document}